\pdfoutput=1
\RequirePackage{ifpdf}
\ifpdf 
\documentclass[pdftex]{sigma}
\else
\documentclass{sigma}
\fi

\newtheorem{Theorem}{Theorem}[section]
\newtheorem*{Theorem*}{Theorem}
\newtheorem{Corollary}[Theorem]{Corollary}
\newtheorem{Lemma}[Theorem]{Lemma}
\newtheorem{Proposition}[Theorem]{Proposition}
\newtheorem{Conjecture}[Theorem]{Conjecture}
\theoremstyle{definition}
\newtheorem{Definition}[Theorem]{Definition}

\newtheorem{Example}[Theorem]{Example}
\newtheorem{Remark}[Theorem]{Remark}

\numberwithin{equation}{section}

\usepackage{bbm}
\usepackage{pgfplots}
\pgfplotsset{compat=1.18}

\renewcommand{\k}{\mathbbm{k}}
\newcommand{\F}{\mathbf{F}}

\newcommand{\sgn}{\mathrm{sign}}

\newcommand{\C}{\mathbb{C}}
\newcommand{\Z}{\mathbb{Z}}

\newcommand{\triv}{\mathrm{triv}}
\newcommand{\std}{\mathrm{std}}
\newcommand{\sign}{\mathrm{sign}}
\newcommand{\trix}{\mathrm{triv-sign-triv}}
\newcommand{\six}{\mathrm{sign-triv-sign}}
\newcommand{\siv}{\mathrm{sign-triv}}
\newcommand{\tign}{\mathrm{triv-sign}}
\newcommand{\trtr}{\mathrm{triv-triv}}

\begin{document}

\newcommand{\arXivNumber}{2412.20673}

\renewcommand{\PaperNumber}{057}

\FirstPageHeading

\ShortArticleName{Hilbert Series of $S_3$-Quasi-Invariant Polynomials in Characteristics 2, 3}

\ArticleName{Hilbert Series of $\boldsymbol{S_3}$-Quasi-Invariant Polynomials\\ in Characteristics 2, 3}

\Author{Frank WANG and Eric YEE}

\AuthorNameForHeading{F.~Wang and E.~Yee}

\Address{Department of Mathematics, Massachusetts Institute of Technology,\\ Cambridge, MA 02139, USA}
\Email{\href{mailto:fw0@mit.edu}{fw0@mit.edu}, \href{mailto:ericyee@mit.edu}{ericyee@mit.edu}}

\ArticleDates{Received January 31, 2025, in final form July 02, 2025; Published online July 13, 2025}

\Abstract{We compute the Hilbert series of the space of $n=3$ variable quasi-invariant polynomials in characteristic $2$ and $3$, capturing the dimension of the homogeneous components of the space, and explicitly describe the generators in the characteristic~2 case. In~doing so we extend the work of the first author in 2023 on quasi-invariant polynomials in characteristic $p>n$ and prove that a sufficient condition found by Ren--Xu in 2020 on when the Hilbert series differs between characteristic $0$ and $p$ is also necessary for $n=3$, $p=2,3$. This is the first description of quasi-invariant polynomials in the case, where the space forms a modular representation over the symmetric group, bringing us closer to describing the quasi-invariant polynomials in all characteristics and numbers of variables.}


\Classification{81R12; 13A50; 20C08}

\section{Introduction}

Let $\k$ be a field, and consider the action of the symmetric group $S_n$ on the space $\k[x_1,\dots,x_n]$ of $\k$-valued polynomials by permuting the variables. A polynomial in $\k[x_1,\dots,x_n]$ is \emph{symmetric} if it is invariant under this action. Equivalently, since $S_n$ is generated by transpositions, a~polynomial~$K$ is symmetric if $s_{i_1i_2}K=K$ or $(1-s_{i_1i_2})K=0$ for all $s_{i_1i_2}\in S_n$. One may consider generalizations of symmetric polynomials in which this condition is relaxed, so that we only require $(1-s_{i_1i_2})K$ be divisible by some large polynomial. This leads to the notion of \emph{quasi-invariant polynomials}.

\begin{Definition}
 Let $\k$ be a field. For $m\in\Z_{\geq 0}$, $n\in\Z_{>0}$, a polynomial $K\in\k[x_1,\dots,x_n]$ is \textit{$m$-quasi-invariant} if for all $s_{i_1i_2}\in S_n$ we have that $(x_{i_1}-x_{i_2})^{2m+1}$ divides $(1-s_{i_1i_2})K$. We denote the space of $m$-quasi-invariants by $Q_m(n,\k)$.
\end{Definition}

Note that the symmetric polynomials are exactly the polynomials that are $m$-quasi-invariant for all $m$. For brevity, we also refer to quasi-invariant polynomials as simply quasi-invariants.

Quasi-invariant polynomials were first introduced by Chalykh and Veselov in 1990 \cite{cmp/1104179957} to describe the harmonic, zero eigenvalue eigenfunctions of quantum Calogero--Moser systems. Calogero--Moser systems are a collection of one-dimensional dynamical particle systems that were found to be both solvable \cite{calogero} and integrable \cite{MOSER1975197}. Due to these properties, they have become extensively studied in mathematical physics, with connections to a number of other fields of mathematics, including representation theory.

Quasi-invariant polynomials were also later found to describe the representation theory of the spherical subalgebra of the rational Cherednik algebra \cite{beg}. This subalgebra is Morita equivalent to the entire rational Cherednik algebra \cite{Etingof2000SymplecticRA}, so quasi-invariants describe representations of rational Cherednik algebras as well. Such algebras have connections to combinatorics, mathematical physics, algebraic geometry, algebraic topology, and more, leading them to become a central topic in representation theory.

Due to these applications, the quasi-invariant polynomials have been studied extensively in recent years. Of particular interest are properties such as its freeness as a module over the symmetric polynomials and the degrees of its generators. To describe these properties, it is useful to consider the Hilbert series of the quasi-invariants, which encapsulates much of this information.

\begin{Definition}
 Let \smash{$V=\bigoplus_{d=0}^\infty V_d$} be a graded vector space. The \textit{Hilbert series} of $V$ is the formal power series
 \[\mathcal{H}(V):=\sum_{d=0}^\infty\dim(V_d)t^d.\]
\end{Definition}

In 2003, Felder and Veselov found the Hilbert series of the space of quasi-invariants in characteristic zero \cite{felder2001action}, proving its freeness in the process. Work on quasi-invariants in characteristic $p$ started in 2020, when Ren and Xu proved a sufficient condition for the Hilbert series of $Q_m(n,\mathbf{F}_p)$ to be different from the Hilbert series of $Q_m(n,\mathbf{Q})$ \cite{Ren_2020}. They accomplished this by computing non-symmetric polynomial ``counterexamples'' in characteristic $p$, where the polynomial has lower degree than any non-symmetric quasi-invariant polynomial in characteristic 0. They also made several conjectures about quasi-invariants in characteristic $p$, including that the condition they found is also sufficient, the quasi-invariants are free, and that the Hilbert polynomial is palindromic for $p>2$. In 2023, the first author proved a general form for the Hilbert series of the quasi-invariants for $n=3$, $p>3$, proving freeness and the palindromicity of the Hilbert polynomial in the process \cite{wang2022explicit}.

We expect the development of the theory of quasi-invariants in characteristic $p$ to be useful in mathematical physics and integrable systems through the theory of $q$-deformed quasi-invariants. These are certain deformations of quasi-invariants in characteristic zero introduced by Chalykh in 2002 \cite{CHALYKH2002193} used to describe eigenfunctions of Macdonald difference operators, which are a generalization of elliptic Calogero--Moser systems \cite{ruij}. We expect the theory of quasi-invariants in characteristic $p$ to be related to the theory of $q$-deformed quasi-invariants when $q$ is a root of unity, in analog to the classical connection between representations of Lie algebras in characteristic $p$ and quantized enveloping algebras \cite{ajs}. We note that a few similarities between these two spaces of quasi-invariants have already been found in \cite{wang2022explicit}.

In this paper, we consider the cases $n=3$, $p=2,3$. These cases differ from the $p>3$ case studied in \cite{wang2022explicit} since in $p=2,3$ the representations of $S_3$ are \emph{modular}, i.e., are not completely reducible. Despite these limitations, we describe the Hilbert series explicitly for all $m$, proving the following.

\begin{Theorem}\label{thm:1}
 Let $\k$ be either $\mathbf{F}_2$ or $\mathbf{F}_3$. Then the Hilbert series for $Q_m(3,\k)$ is given by
 \[\mathcal{H}(Q_m(3,\k))=\frac{1+2t^d+2t^{6m+3-d}+t^{6m+3}}{(1-t)(1-t^2)(1-t^3)},\]
 where $d=3m+1$ if there is no Ren--Xu counterexample and $d$ is the degree of the minimal degree Ren--Xu counterexample otherwise. In particular, the conditions found in {\rm \cite{Ren_2020}} for the Hilbert series of $Q_m(3,\k)$ to be different from the Hilbert series of $Q_m(3,\mathbf{Q})$ are necessary.
\end{Theorem}

Note that this result also implies freeness and the palindromicity of the Hilbert polynomial.

In the case $p=2$, we also define $m$-quasi-invariants in the case where $m$ is a half-integer and prove an analogous statement to Theorem~\ref{thm:1} in this case. Using quasi-invariants at half-integers, we also compute the generators of $Q_m(3,\mathbf{F}_2)$ as an $\mathbf{F}_2[x_1,x_2,x_3]^{S_3}$-module explicitly.

In Section \ref{prelim}, we state some of the basic facts about quasi-invariant polynomials and introduce modular representations of $S_3$. In Section \ref{sec:p=2}, we compute the generators of $Q_m(3,\mathbf{F}_2)$, proving Theorem~\ref{thm:1} for $p=2$ in the process. In Section \ref{sec:wang}, we begin discussing $p=3$, and show that some properties of quasi-invariants in 3 variables from \cite{wang2022explicit} carry over to the $p=3$ case after converting from the standard representation to the $\siv$ representation. In Section~\ref{sec:renxu}, we show that minimal degree Ren--Xu counterexamples are the lowest degree non-symmetric generators for $Q_m(3,\mathbf{F}_3)$ and show that there is one other higher degree generator belonging to the $\siv$ representation. Finally, in Section \ref{sec:reps}, we consider all other indecomposable representations of $S_3$ in $Q_m(3,\mathbf{F}_3)$, finishing the proof of Theorem~\ref{thm:1} for $p=3$.

\section{Preliminaries}\label{prelim}
We start with some useful properties of the quasi-invariants.
\begin{Proposition}[\cite{etingof2002lectures}]\label{firstProp}
 Let $\k$ be a field.
\begin{enumerate}\itemsep=0pt
 \item[$1$.] $\k[x_1,x_2,x_3]^{S_3}\subset Q_m(3,\k)$, $Q_0(3,\k) = \k[x_1,x_2,x_3]$, and $Q_m(3,\k)\supset Q_{m'}(3,\k)$, where ${m'>m}$.
 \item[$2$.] $Q_m(3,\k)$ is a ring.
 \item[$3$.] $Q_m(3,\k)$ is a finitely generated $\k[x_1,x_2,x_3]^{S_3}$-module.
\end{enumerate}
\end{Proposition}

Note that \cite{etingof2002lectures} proves Proposition~\ref{firstProp} in the case, where $\k=\C$. However, the proofs for the first two assertions work over any field, and the last assertion follows from the Hilbert basis theorem. In view of the structure of $Q_m(3,\k)$ as a module over the symmetric polynomials, given some $K\in Q_m(3,\k)$, we will frequently refer to quasi-invariant polynomials that can be obtained via scalar multiplication of $Q$ by a symmetric polynomial. To distinguish these polynomials from the ordinary $\k$-multiples of $Q$, we will refer to them as \emph{symmetric polynomial multiples} of $Q$.

We consider $Q_m(3, \mathbf{F}_{2})$ and $Q_m(3, \mathbf{F}_{3})$ as representations of $S_3$, where $S_3$ permutes the variables $x_1$, $x_2$, $x_3$. Since $Q_m(3, \mathbf{F}_{2})$ and $Q_m(3, \mathbf{F}_{3})$ are vector spaces over $\mathbf{F}_2$ and $\mathbf{F}_3$ respectively and the characteristics $2$ and $3$ divide $|S_3|$, $Q_m(3, \mathbf{F}_{2})$ and $Q_m(3, \mathbf{F}_{3})$ are modular representations~of~$S_3$.

\begin{Proposition}
 $Q_m(3, \mathbf{F}_{2})$ and $Q_m(3, \mathbf{F}_{3})$ are modular representations of $S_3$.
\end{Proposition}

First, we consider characteristic $2$.

\subsection[Preliminary definitions for p=2]{Preliminary definitions for $\boldsymbol{p=2}$}

We describe the indecomposable and irreducible representations of $S_3$ for $p=2$.

\begin{Proposition}[\cite{Alperin_1986}]\label{prop:p2indecomp}
 There are $3$ irreducible or indecomposable representations of $S_{3}$ in characteristic $2$:
 \begin{enumerate}\itemsep=0pt
 \item[$1$.] $\triv$ is the irreducible representation of $S_3$ that is acted on trivially by $S_{3}$.
 \item[$2$.] $\std$ is the $2$-dimensional irreducible representation of $S_3$ obtained by reducing the standard representation in characteristic $0$ mod~$2$.
 \item[$3$.] $\trtr$ is the $2$-dimensional indecomposable representation that contains a copy of $\triv$ as a subrepresentation such that the quotient of $\trtr$ by this subrepresentation is $\triv$.
 \end{enumerate}
\end{Proposition}

\begin{Example}\label{eg:p2reps}
 The polynomial $E_{\trtr}:=x_1^2x_2+x_2^2x_3+x_3^2x_1\in\mathbf{F}_2[x_1,x_2,x_3]$ generates a~copy of $\trtr$. To see this, note that for any $i_1$, $i_2$, we have
 \[(1-s_{i_1i_2})E_\trtr=x_1^2x_2+x_1x_2^2+x_1^2x_3+x_1x_3^2+x_2^2x_3+x_2x_3^2\in\mathbf{F}_2[x_1,x_2,x_3]^{S_3}.\]
 Since the transpositions generate $S_3$, $E_\trtr$ generates a two-dimensional representation that contains $\triv$ as a subrepresentation. Moreover, since $E_\trtr$ is not symmetric, this representation is not $\triv\oplus\triv$, so it must be $\trtr$.
\end{Example}

We then study the behaviors of each indecomposable representation in the quasi-invariants. We define $Q_m(3,\mathbf{F}_2)_{\triv}$ and $Q_m(3,\mathbf{F}_2)_{\std}$ to be the direct sum of all copies of $\triv$ and $\std$ respectively in the quasi-invariants. We also define $Q_m(3,\mathbf{F}_2)_\trtr$ to be the direct sum of all copies of $\triv$ and $\trtr$.

\begin{Remark}
 We cannot define $Q_m(3,\mathbf{F}_2)_\trtr$ to exclude copies of $\triv$ since we can add elements of $Q_m(3,\mathbf{F}_2)_\triv$ to copies of $\trtr$ and still obtain a copy of $\trtr$. For example, \smash{$F:=E_\trtr+x_1^3+x_2^3+x_3^3$} still satisfies $(1-s_{i_1i_2})F=(1-s_{i_1i_2})E_\trtr$ for all~$i_1$,~$i_2$, so it generates a copy of $\trtr$ by the same argument as Example~\ref{eg:p2reps}.
\end{Remark}

\begin{Proposition}[{\cite{wang2022explicit}}]\label{prop:p2triv}
 As an $\mathbf{F}_2[x_1,x_2,x_3]^{S_3}$-module, $Q_m(3,\mathbf{F}_2)_{\triv}$ is freely generated by $1$.
\end{Proposition}

Note that by the classification of indecomposables in Proposition~\ref{prop:p2indecomp}, every extension of $\std$ and every extension of a module by $\std$ splits. Thus $Q_m(3,\mathbf{F}_2)_\std$ is a direct summand of $Q_m(3,\mathbf{F}_2)$ (whose complement is $Q_m(3,\mathbf{F}_2)_\trtr$), and we mainly consider $Q_m(3,\mathbf{F}_2)_\std$. $Q_m(3,\mathbf{F}_2)_{\std}$ is generated as a $\mathbf{F}_2[x_1,x_2,x_3]^{S_3}$-module by homogeneous copies of $\std$, so following~\cite{wang2022explicit}, we consider \textit{generating representations} of $Q_m(3,\mathbf{F}_2)_\std$ as homogeneous copies of $\std$ in a generators and relations presentation of $Q_m(3,\mathbf{F}_2)_{\std}$ with a minimal generator set.

\subsubsection{Quasi-invariants at half-integers}

Note that if $\k$ is a field with $\mathrm{char\,}\k\neq 2$ and $m\in\Z_{\geq 0}$, then for any $K\in\k[x_1,\dots,x_n]$, $(x_{i_1}-x_{i_2})^{2m}|(1-s_{i_1i_2})K$ implies $(x_{i_1}-x_{i_2})^{2m+1}|(1-s_{i_1i_2})K$ since $(1-s_{i_1i_2})K$ is $s_{i_1i_2}$-antiinvariant, hence the exponent $2m+1$ in the definition of quasi-invariant polynomials. But this does not hold in characteristic 2, since there is no concept of antiinvariants. Indeed, one can check that for $K=x_1^2+x_2^2$, we have $(x_{i_1}-x_{i_2})^2|(1-s_{i_1i_2})K$ for all $i_1$, $i_2$, but $(x_{i_1}-x_{i_2})^3\nmid|(1-s_{i_1i_2})K$ if~$i_1=1,2$,~$i_2\neq 1,2$.

We encapsulate this data by extending the definition of quasi-invariants to half-integers when ${p=2}$. For example, $K=x_1^2+x_2^2$ is \smash{$\frac{1}{2}$}-quasi-invariant, and this is in fact the minimal degree nonsymmetric \smash{$\frac{1}{2}$}-quasi-invariant polynomial. Proposition~\ref{firstProp} still holds when $m,m'$ are half-integers, and the definitions of $Q_m(3,\mathbf{F}_2)_\triv$, $Q_m(3,\mathbf{F}_2)_\std$ also naturally extend to half-integer~$m$. So from now on, whenever we refer to quasi-invariants in characteristic $2$ we let $m$ be a~half-integer.

\subsection[Preliminary definitions for p=3]{Preliminary definitions for $\boldsymbol{p=3}$}

Next, we define the indecomposable and irreducible representations of $S_3$.

\begin{Proposition}[\cite{Alperin_1986}]
 There are $6$ indecomposable or irreducible representations in $S_3$ in characteristic $3$:
 \begin{enumerate}\itemsep=0pt
 \item[$1$.] $\triv$ is the irreducible representation of $S_3$ that is acted on trivially by $S_3$.
 \item[$2$.] $\sign$ is the irreducible representation of $S_3$ that is acted on by negation by the transpositions.
 \item[$3$.] $\siv$ is the indecomposable representation that contains a copy of $\triv$ as a subrepresentation, such that the quotient of $\siv$ by this subrepresentation is $\sign$.
 \item[$4$.] $\tign$ is the indecomposable representation that contains a copy of $\sign$ as a subrepresentation, such that the quotient of $\tign$ by this subrepresentation is $\triv$.
 \item[$5$.] $\trix$ is the indecomposable representation that contains a copy of $\siv$ as a subrepresentation, such that the quotient of $\trix$ by this subrepresentation is $\triv$.
 \item[$6$.] $\six$ is the indecomposable representation that contains a copy of $\tign$ as a subrepresentation, such that the quotient of $\six$ by this subrepresentation is $\sign$.
 \end{enumerate}
\end{Proposition}

Provided are some examples of copies of these indecomposable representations:

\begin{Example}\label{exSiv}
 The space $W\subset\F_3[x_1,x_2,x_3]$ spanned by $x_1+x_2+x_3$ and $x_1-x_2$ over $\mathbf{F}_3$ is copy of $\siv$. Indeed, the space $T\subset W$ spanned by $x_1+x_2+x_3$ is a copy of $\triv$. One can check $x_1-x_2\in W/T$ is acted by negation by all transpositions in $S_3$ and $W/T$ is $1$-dimensional so $W/T$ is a copy of $\sign$. Finally, it is easy to show that there are no copies of $\triv$ or $\sign$ in $W$ other than $T$. Since $V$ has a unique irreducible subrepresentation, it is indecomposable, and we conclude that it is a copy of $\siv$.
\end{Example}

\begin{Example}\label{deg1polys}
 The space $V\subset\F_3[x_1,x_2,x_3]$ consisting of homogeneous linear polynomials is a~copy of $\trix$. Indeed, $W\subset V$ from Example~\ref{exSiv} is a copy of $\siv$. Then~$V/W$ is one-dimensional, and one can check that it is a copy of $\triv$. Finally, it is easy to show that there are no copies of $\triv$ or $\sign$ in $V$ other than $T$, so $V$ has a unique irreducible subrepresentation, it is indecomposable, and we conclude that it is a copy of $\trix$.
\end{Example}

\begin{Example}
 Similarly, one may check that the space $U$ spanned by
 \[(x_1-x_2)(x_1-x_3)(x_2-x_3),\qquad-x_1^2x_2-x_1^2x_3+x_1x_2^2+x_1x_3^2\]
 over $\mathbf{F}_3$ is a copy of $\tign$ and that the space spanned by
 \[(x_1-x_2)(x_1-x_3)(x_2-x_3),\qquad-x_1^2x_2-x_1^2x_3+x_1x_2^2+x_1x_3^2,\qquad(x_1-x_2)x_1x_2\]
 is a copy of $\six$.
\end{Example}

Similarly to the $p=2$ case, we define $Q_m(3,\mathbf{F}_3)_{\sign}$ and $Q_m(3,\mathbf{F}_3)_{\triv}$ to be the direct sum of all copies of $\sign$ and $\triv$ in $Q_m(3,\mathbf{F}_3)$, respectively.
\begin{Proposition}[\cite{wang2022explicit}]\label{prop:trivsign}
 As $\mathbf{F}_3[x_1,x_2,x_3]^{S_3}$-modules,
 \begin{enumerate}\itemsep=0pt
 \item[$1$.] $Q_m(3,\mathbf{F}_3)_\triv$ is freely generated by $1$.
 \item[$2$.] $Q_m(3,\mathbf{F}_3)_\sgn$ is freely generated by $\prod_{i_1<i_2}(x_{i_1}-x_{i_2})^{2m+1}$.
 \end{enumerate}
\end{Proposition}

Next we define $Q_m(3,\mathbf{F}_3)_{\siv}$ as the direct sum of all copies of $\sign$, $\triv$, and $\siv$. For this paper we consider generators of $Q_m(3,\mathbf{F}_3)_{\siv}$ to be homogeneous polynomials other than $1$ and \smash{$\prod_{i_1<i_2}(x_{i_1}-x_{i_2})^{2m+1}$} such that they are in the $(-1)$-eigenspace of $s_{12}$ and are in a~generators and relations presentation of $Q_m(3,\mathbf{F}_3)_{\siv}$ as an $\mathbf{F}_3[x_1,x_2,x_3]^{S_3}$-module with the least number of generators. Moreover, if $K$ is a generator of $Q_m(3,\mathbf{F}_3)_{\siv}$ then it necessarily generates a copy of $\siv$ since we assumed $K$ neither generates $\triv$ nor $\sign$.

\begin{Remark}
 Similar to in the $p=2$ case, we cannot define $Q_m(3,\mathbf{F}_3)_{\siv}$ to exclude copies of $\sign$ since we can add elements of $Q_m(3,\mathbf{F}_3)_{\sign}$ to copies of $\siv$ and still obtain a copy of $\siv$. For example, the spaces spanned by
 \[\big(x_1^6-x_2^6\big)(x_1+x_2+x_3)^3,\qquad\big(x_1^6+x_2^6+x_3^6\big)(x_1+x_2+x_3)^3\]
 and
 \[\prod_{i_1<i_2}(x_{i_1}-x_{i_2})^{3}+\big(x_1^6-x_2^6\big)(x_1+x_2+x_3)^3,\qquad\big(x_1^6+x_2^6+x_3^6\big)(x_1+x_2+x_3)^3\]
 generate two copies of $\siv$ in $Q_1(3,\F_3)$, and their sum contains
 \[\prod_{i_1<i_2}(x_{i_1}-x_{i_2})^{3}\in Q_1(3,\mathbf{F}_3)_{\sign}.\]
\end{Remark}

\begin{Remark}
 One could define subspaces of $Q_m(3,\mathbf{F}_3)$ for $\trix$, $\six$, $\tign$ similar to $Q_m(3,\mathbf{F}_3)_{\siv}$, however this is not particularly helpful, as unlike for ${p=2}$, we cannot decompose $Q_m(3,\mathbf{F}_3)$ into a direct sum of subspaces of this form. The space $Q_m(3,\mathbf{F}_3)_{\siv}$ is still relevant, as it is the critical piece to understanding quasi-invariants in characteristic 3, as we see in Sections \ref{sec:wang} and \ref{sec:renxu}.
\end{Remark}

\section[Quasi-invariants in characteristic 2]{Quasi-invariants in characteristic 2}\label{sec:p=2}

In this section, we write down explicit generators for $Q_m(3,\mathbf{F}_2)$ and prove Theorem~\ref{thm:1} for $p=2$. Note that we already know the structure of $Q_m(3,\mathbf{F}_2)_\triv$ from Proposition~\ref{prop:p2triv}. We start by extending this to $Q_m(3,\mathbf{F}_2)_\trtr$.

\begin{Proposition}\label{prop:trivtriv}
 As an $\mathbf{F}_2[x_1,x_2,x_3]^{S_3}$-module, $Q_m(3,\mathbf{F}_2)_{\triv-\triv}$ is freely generated by $1$ and $E_{\trtr}\prod(x_{i_1}-x_{i_2})^{2m}$.
\end{Proposition}
\begin{proof}
 Let $K$ be a nonsymmetric element of $Q_m(3,\mathbf{F}_2)_{\triv-\triv}$ so that $(x_{i_1}-x_{i_2})^{2m+1}$ divides $(1+s_{i_1i_2})K$. Because
 \[(1+s_{12})K=(1+s_{13})K=(1+s_{23})K,\]
 we have $(1+s_{i_1i_2})K = P\prod (x_{i_1}-x_{i_2})^{2m+1}$ for some symmetric polynomial $P$. Letting $G=E_{\trtr}\prod (x_{i_1}-x_{i_2})^{2m}$ yields $(1+s_{i_1i_2})G = \prod (x_{i_1}-x_{i_2})^{2m+1}$. Thus $(1+s_{i_1i_2})PG = (1+s_{i_1i_2})K$ and $(1+s_{i_1i_2})(PG-K)=0$, so $PG-K$ is symmetric and $K$ is generated by $G$ and $1$. Moreover, since $G$ is not symmetric, $P$ and $G$ have no relation implying freeness.
\end{proof}

We have an explicit description of $Q_m(3,\mathbf{F}_2)_\trtr$, so it remains to compute the generators and relations of $Q_m(3,\mathbf{F}_2)_\std$. A number of the properties of $Q_m(3,\mathbf{F}_p)$ for $p>3$ found in \cite{wang2022explicit} are true for $Q_m(3,\mathbf{F}_2)$. We prove these first.

If $V$ is a copy of $\std$, then we denote by $V_{i_1i_2}$ the $1$-eigenspace of $s_{i_1i_2}$ in $V$.

\begin{Lemma}
 Let $V$ be a copy of $\std$ in $Q_m(3,\mathbf{F}_2)_{\std}$, and let $K\in V_{i_1i_2}$. Then we have ${K+sK+s^2K}=0$, where $s=(1\,2\,3)\in S_3$ and $K=(x_{i_1}-x_{i_2})^{2m+1}K'$ for some polynomial~$K'$ that is invariant under the action of $s_{i_1i_2}$. Conversely, let $K'$ be an $s_{12}$-invariant polynomial such that
 \[(x_1-x_2)^{2m+1}K'+(x_2-x_3)^{2m+1}sK'+(x_3-x_1)^{2m+1}s^2K'=0.\]
 Then $(x_1-x_2)^{2m+1}K'$ belongs to the $1$-eigenspace of $s_{12}$ in some copy of $\std$ inside $Q_m(3,\mathbf{F}_2)_\std$.
\end{Lemma}
\begin{proof}
 For the first statement, $K+sK+s^2K=0$ holds for any copy of $\std$. For the next, suppose $\{i_1,i_2,i_3\} = \{1,2,3\}$ for some integer $i_3$. Then $(1-s_{i_1i_3})K = s_{i_2i_3}K$, so $(x_{i_1}-x_{i_3})^{2m+1}|s_{i_2i_3}K$, implying $(x_{i_1}-x_{i_2})^{2m+1}|K$.
 The second statement follows from the proof in~\cite{wang2022explicit}.
\end{proof}

\begin{Corollary}
 Let $V$ be a generating representation of $Q_m(3,\k)_\std$ and let $K\in V_{i_1i_2}$. Let us write $K=(x_{i_1}-x_{i_2})^{2m+1}K'$. Then $K'$ is not divisible by any nonconstant symmetric polynomial.
\end{Corollary}

The proof of this statement is identical to the one in \cite{wang2022explicit}.

\begin{Lemma}\label{char3.4}
 Let $V$, $W$ be distinct generating representations of $Q_m(3,\k)_\std$. Let ${K\in V_{12}}$, ${L\in W_{12}}$. For $\sigma K\sigma L:=(\sigma K)(\sigma L)$, we have that $KL+s_{13}Ks_{23}L$ is a nonsymmetric element of $Q_m(3,\k)_\trtr$ and $\deg V+\deg W\geq 6m+3$.
\end{Lemma}
\begin{proof}
 $KL+s_{13}Ks_{23}L$ is an element of $Q_m(3,\mathbf{F}_2)$ since the quasi-invariants form a ring by Proposition~\ref{firstProp}. Using that $s_{12}K=K$ and $s_{12}L=L$, we have that
 \begin{gather*}
 (1+s_{12})(KL+s_{13}Ks_{23}L) = s_{23}Ks_{13}L+s_{13}Ks_{23}L,\\
 (1+s_{13})(KL+s_{13}Ks_{23}L) = KL+s_{13}Ks_{23}L+s_{13}Ks_{13}L + Ks_{23}L = Ks_{13}L+s_{13}KL,\\
 (1+s_{23})(KL+s_{13}Ks_{23}L) = KL+s_{13}Ks_{23}L+s_{23}Ks_{23}L + s_{13}KL = Ks_{23}L+s_{23}KL.
 \end{gather*}
 One can check that each polynomial is a transposition of another and that they are symmetric due to the structure of $\trtr$, so they are all the same symmetric polynomial. Thus ${KL+s_{13}Ks_{23}L}$ lies in a quotient of a copy of $\trtr$. Note that by the same argument as in \cite{wang2022explicit}, we have $Ks_{23}L+s_{23}KL\neq 0$, so $KL+s_{13}Ks_{23}L$ is nonsymmetric and must generate a~copy of $\trtr$.

 By Proposition~\ref{prop:trivtriv}, $KL+s_{13}Ks_{23}L$ has degree at least $6m+3$, so $\deg V+\deg W\geq 6m+3$ as desired.
\end{proof}

\begin{Lemma}\label{charFree}
 Assume that there exist generating representations $V$, $W$ of $Q_m(3,\mathbf{F}_2)_\std$ such that $\deg V+\deg W=6m+3$. Then $Q_m(3,\mathbf{F}_2)_\std$ is a free module over $\k[x_1,x_2,x_3]^{S_3}$ generated by~$V$ and $W$.
\end{Lemma}
\begin{proof}
 Assume for the sake of contradiction there exists another generator $U$ of $Q_m(3,\mathbf{F}_2)_{\std}$. Supposing $\deg W\geq \deg V$, by Lemma~\ref{char3.4}, $\deg U\geq \deg W$. By Lemma~\ref{char3.4}, if $K\in V_{12}$, $L\in W_{12}$, and $T\in U_{12}$ then $KL+s_{13}Ks_{23}L$ and $KT+s_{13}Ks_{23}T$ are both nonsymmetric elements of $Q_m(3,\mathbf{F}_2)_{\trtr}$. Moreover, we have
 \[(1+s_{12})(KL+s_{13}Ks_{23}L) = s_{23}Ks_{13}L+s_{13}Ks_{23}L = \prod (x_{i_1}-x_{i_2})^{2m+1},\]
 and
 \[(1+s_{12})(KT+s_{13}Ks_{23}T) = s_{23}Ks_{13}T+s_{13}Ks_{23}T = Q\prod (x_{i_1}-x_{i_2})^{2m+1}\]
 for some symmetric polynomial $Q$. From there we may proceed identically to \cite{wang2022explicit}.
\end{proof}

Now, we are ready to prove Theorem~\ref{thm:1} for $p=2$.

\begin{Theorem}\label{thm:p2}
 Let $a$ be the largest natural number such that $2^a< 2m+1$. Then $Q_m(3,\mathbf{F}_2)_{\std}$ is freely generated by $(x_1-x_2)^{2^{a+1}}$ and $(x_1-x_2)^{2^a}\prod (x_{i_1}-x_{i_2})^{2m+1-2^a}$.
\end{Theorem}

\begin{Remark}
 Note that when $m$ is an integer, the degrees of the generators in this theorem agree with the degrees conjectured in \cite{Ren_2020}. In particular, when $2^{a+1}$ is one of $3m+1$, ${3m+2}$, we actually have that the Hilbert series of $Q_m(3,\F_2)$ and $Q_m(3,\mathbf{Q})$ agree, so $(x_1-x_2)^{2^{a+1}}$, \smash{$(x_1-x_2)^{2^a}\prod (x_{i_1}-x_{i_2})^{2m+1-2^a}$} are the reductions modulo 2 of the generators of $Q_m(3,\mathbf{Q})$, when written as integer polynomials with coprime coefficients.
\end{Remark}

\begin{proof}[Proof of Theorem~\ref{thm:p2}]
 We prove this by induction on $m$.

 The generators of $Q_0(3,\mathbf{F}_2)_{\std}$ are $(x_1-x_2)$ and $(x_1-x_2)^2$, completing our base case.

 Let $j$ be a half-integer, and suppose that \smash{$Q_{j-\frac{1}{2}}(3,\mathbf{F}_2)_\std$} is freely generated by $(x_1-x_2)^{2^{a+1}}$ and \smash{$(x_1-x_2)^{2^a}\prod (x_{i_1}-x_{i_2})^{2j-2^a}$},
 where $2^a$ is the greatest such power of $2$ less than $2j$. If~${2j\neq 2^{a+1}}$, then $2^{a}$ is the largest power of $2$ less than $2j+1$, so $(x_1-x_2)^{2^{a+1}}$ and $(x_1-x_2)^{2^a}\prod (x_{i_1}-x_{i_2})^{2j+1-2^a}$ are both in $Q_j(3,\mathbf{F}_2)$. Further, $(x_1-x_2)^{2^{a+1}}$ must be a generator and if $(x_1-x_2)^{2^a}\prod (x_{i_1}-x_{i_2})^{2j+1-2^a}$ is a not a generator, by Lemma~\ref{char3.4}, $(x_1-x_2)^{2^a}\prod (x_{i_1}-x_{i_2})^{2j+1-2^a}$ is generated by $(x_1-x_2)^{2^{a+1}}$ which implies a relation between $(x_1-x_2)^{2^{a+1}}$ and ${(x_1-x_2)^{2^a}\prod (x_{i_1}-x_{i_2})^{2j-2^a}}$. Because they freely generate \smash{$Q_{j-\frac{1}{2}}(3,\mathbf{F}_2)$}, this is impossible. Thus $(x_1-x_2)^{2^{a+1}}$ and $(x_1-x_2)^{2^a}\smash{\prod (x_{i_1}-x_{i_2})^{2j+1-2^a}}$ freely generate $Q_j(3,\mathbf{F}_2)_{\std}$ by Lemma~\ref{charFree}.

 If $2j= 2^{a+1}$, then both $(x_2-x_3)^{2^{a+1}}$ and $(x_2-x_3)^{2^{a+2}}\prod (x_{i_1}-x_{i_2})^{2j+1-2^a}$ lie in $Q_j(3,\mathbf{F}_2)_{\std}$. The former is a generator by our inductive hypothesis.
 Since $2^{a+1}+2^{a+2}+3=6j+3$, if the latter is not a generator, then by Lemma~\ref{char3.4}, $(x_2-x_3)^{2^{a+2}\prod (x_{i_1}-x_{i_2})^{2j+1-2^a}}$ is generated by $(x_2-x_3)^{2^{a+1}}$, which is false. Thus $(x_2-x_3)^{2^{a+1}}$ and $(x_2-x_3)^{2^{a+2}}\prod (x_{i_1}-x_{i_2})^{2j+1-2^a}$ freely generate $Q_j(3,\mathbf{F}_2)_{\std}$ by Lemma~\ref{charFree} as desired.
\end{proof}

\section[Properties of 3 variable quasi-invariants]{Properties of $\boldsymbol{3}$ variable quasi-invariants}\label{sec:wang}
Similarly to the $p=2$ case, we can adapt many of the properties of $Q_m(3,\mathbf{F}_p)$ for $p>3$ found in \cite{wang2022explicit} to the $p=3$ case. We accomplish this by converting $\std$ to $\siv$. For example, in~$Q_0(3,\mathbf{F}_p)$ for $p>3$, the space spanned by $x_1-x_2$, $x_1-x_3$ is a copy of $\std$. However, in~$Q_0(3,\mathbf{F}_3)$, the space spanned by $x_1-x_2$, $x_1-x_3$ becomes a copy of $\siv$. Using this, we may show that there are equivalents of Lemmas 3.2--3.5 from \cite{wang2022explicit} in characteristic~$3$.

We define $V_{i_1i_2}^{-}$ to be the $(-1)$-eigenspace of $s_{i_1i_2}$ in $V$, where $V$ is a copy of $\std$ or $\siv$. Note that if $v\in V_{i_1i_2}^{-}$ we have $v=s_{23}v+s_{13}v$. The following lemma and corollary correspond to Lemma~3.2 and Corollary~3.3 from \cite{wang2022explicit}, respectively.

\begin{Lemma}\label{lem:3.1}
 Let $V$ be a copy of $\siv$ in $Q_m(3,\mathbf{F}_3)_{\siv}$, and let $K\in V_{i_1i_2}^{-}$. Then we have $K+sK+s^2K=0$, where $s=(1\,2\,3)\in S_3$ and $K=(x_{i_1}-x_{i_2})^{2m+1}K'$ for some polynomial~$K'$ that is invariant under the action of $s_{i_1i_2}$. Conversely, let $K'$ be an $s_{12}$-invariant polynomial such that
 \[(x_1-x_2)^{2m+1}K'+(x_2-x_3)^{2m+1}sK'+(x_3-x_1)^{2m+1}s^2K'=0.\]
 Then $(x_1-x_2)^{2m+1}K'$ either belongs to $Q_m(3,\mathbf{F}_3)_\sign$ or the $(-1)$-eigenspace of $s_{12}$ in some copy of $\siv$ inside $Q_m(3,\mathbf{F}_3)_\siv$.
\end{Lemma}

\begin{proof}
 The proof is largely the same as in \cite{wang2022explicit}; the only difference is in the last step. Namely, now we have 2 2-dimensional indecomposable representations $\siv$ and $\tign$, but an~element in the $(-1)$-eigenspace of $s_{12}$ in $\tign$ must be in a copy of $\sign$.
\end{proof}

\begin{Corollary}\label{cor:3.2}
 Let $K$ be a generator of $Q_m(3,\mathbf{F}_3)_\siv$ in $V_{i_1i_2}^-$ for some copy $V$ of $\siv$ and write $K=(x_{i_1}-x_{i_2})^{2m+1}K'$. Then $K'$ is not divisible by any nonconstant symmetric polynomial.
\end{Corollary}

The proof of this corollary is identical to the proof of \cite[Corollary 3.3]{wang2022explicit}.

We define generators of $Q_m(3,\mathbf{F}_3)$ to be ``distinct'' if they are either in different degrees, or if no linear combination of them over $\mathbf{F}_3$ is generated by lower degree generators.

\begin{Lemma}\label{lem:3.3}
 Let $K$ and $L$ be distinct generators of $Q_m(3,\k)_\siv$, and let $V$ and $W$ be the~copies of $\siv$ generated by $K$ and $L$ respectively such that $K\in V_{i_1i_2}^{-}$ and $L\in W_{i_1i_2}^{-}$. Then $Ks_{23}L - Ls_{23}K$ is a nonzero element of $Q_m(3,\mathbf{F}_3)_{\sign}$ and $\deg V+\deg W\geq 6m+3$.
\end{Lemma}

Noting that $\wedge^2(\siv) = \sign$, the proof of this lemma is also identical to the proof of~\mbox{\cite[Lemma 3.4]{wang2022explicit}}.

Lemma 3.5 from \cite{wang2022explicit} does not completely hold in characteristic $3$. A very similar and useful version does, however, and we have the following.

\begin{Lemma}\label{lem:3.4}
 Assume that there exists generators $K$ and $L$ of $Q_m(3,\mathbf{F}_3)_\siv$ such that $\deg K+\deg L=6m+3$. Then $Q_m(3,\mathbf{F}_3)_\siv$ is freely generated by $K$, $L$, and $1$ over $\mathbf{F}_3[x_1,x_2,x_3]^{S_3}$.
\end{Lemma}
\begin{proof}
 We note that $(L+s_{23}L)K-(K+s_{23}K)L=c\prod_{i_1<i_2}(x_{i_1}-x_{i_2})^{2m+1}$ for some $c\neq 0$ by Lemma~\ref{lem:3.3}. Moreover, $L+s_{23}L$ and $K+s_{23}K$ are symmetric because $K$ and $L$ are both acted on by negation by $s_{12}$, so elements in $Q_m(3,\mathbf{F}_3)_{\sign}$ are generated by $K$ and $L$. From there, the fact that $Q_m(3,\mathbf{F}_3)_\siv$ is generated by $K$, $L$, and $1$ over $\mathbf{F}_3[x_1,x_2,x_3]^{S_3}$ follows from the first part of the proof from \cite{wang2022explicit}.

 To prove freeness, assume for the sake of contradiction that there was a relation $PK+QL+S=0$ for symmetric polynomials $P$, $Q$, and $S$. $PK$ and $QL$ are both in the $(-1)$-eigenspace of $s_{12}$ while $S$ is not, so $S=0$. Thus we have $PK=-QL$ and from there we can proceed the same as \cite{wang2022explicit}.

\end{proof}

\section{Ren--Xu counterexamples}\label{sec:renxu}

We aim to explicitly describe the Hilbert series of $Q_m(3,\mathbf{F}_3)$. To do so we wish to identify the generators of $Q_m(3,\mathbf{F}_3)_{\siv}$.

In \cite{Ren_2020}, Ren and Xu found polynomials of the form $P_k^{3^a}\prod (x_{i_1}-x_{i_2})^{2b}$ in $Q_m(3,\mathbf{F}_3)$ with degree strictly less than $3m+1$, where $P_k$ is the map of the $3k+1$ degree generator of $Q_k(3,\mathbf{Q})$ into characteristic $3$ and where $a$, $k$, and $b$ are natural numbers. We refer to these polynomials as Ren--Xu counterexamples as they demonstrate the Hilbert series of $Q_m(3,\mathbf{F}_3)$ differs from that of $Q_m(3,\mathbf{Q})$ for certain $m$.

\begin{Definition}
 Let $\overline{P_k}$ be the generator of $Q_k(3,\mathbf{Q})$ of degree $3k+1$ in the $(-1)$-eigenspace of $s_{12}$, expressed as an integer polynomial with coprime coefficients. Let $P_k$ be the image of $\overline{P_k}$ under the quotient map $\Z[x_1,x_2,x_3]\to\F_3[x_1,x_2,x_3]$. Define the set $X$ as the set of all natural numbers $m$ such that $Q_m(3,\mathbf{F}_3)$ has a Ren--Xu counterexample. Let $R_m$ be a lowest degree Ren--Xu counterexample in $Q_m(3,\mathbf{F}_3)$ for all $m\in X$.
\end{Definition}

A key step in describing the Hilbert series of $Q_m(3,\mathbf{F}_3)$ is proving Ren--Xu's conjecture~\cite{Ren_2020} for ${n=3}$ and~${p=3}$.

\begin{Conjecture}[\cite{Ren_2020}]\label{conj:ren}
 If the Hilbert series of $Q_m(n,\mathbf{F}_p)$ differs from that of $Q_m(n,\mathbf{Q})$, then there exists integers $a\geq 0$ and $k\geq 0$ such that
 \[\frac{mn(n-2)+\binom{n}{2}}{n(n-2)k+\binom{n}{2}-1}\leq p^a\leq \frac{mn}{nk+1}.\]
\end{Conjecture}

The main step for proving the conjecture for $n=3$, $p=3$ is the following theorem.

\begin{Theorem}\label{justRens}
 $Q_m(3,\mathbf{F}_3)_{\siv}$ is either freely generated by a generator of degree $3m+1$, $3m+2$, and the polynomial $1$, or it is freely generated by $R_m$, another generator in degree $6m+3-\deg R_m$, and the polynomial $1$.
\end{Theorem}

To prove this theorem, we first describe the Ren--Xu counterexamples.

\begin{Lemma}\label{nonExsOnly}
 If $m\in X$, we must have $R_m=P_k^{3^a}\prod (x_{i_1}-x_{i_2})^{2b}$, where $a$, $b$, $k$ are natural numbers and $k\not \in X$.
\end{Lemma}
\begin{proof}
 Assume for contradiction that there exists a nonnegative integer $m\in X$ such that $R_m=P_k^{3^a}\prod (x_{i_1}-x_{i_2})^{2b}$, where $a$, $b$, $k$ are natural numbers and $k\in X$. Then if $R_k=P_{l}^{3^c}\prod (x_{i_1}-x_{i_2})^{2d}$, the polynomial
 \[{R_k}^{3^a}\prod (x_{i_1}-x_{i_2})^{2b} = P_{l}^{3^{a+c}}\prod (x_{i_1}-x_{i_2})^{2d\cdot 3^a + 2b}\]
 has a strictly smaller degree than $R_m$ since $\deg R_k<3k+1=\deg P_k$. Moreover, it is at least $m$-quasi-invariant, so it is a Ren--Xu counterexample for $Q_m(3,\mathbf{F}_3)$. Yet $R_m$ is a minimal counterexample, giving a contradiction.
\end{proof}

This lemma allows us to consider only counterexamples $P_k^{3^a}\prod (x_{i_1}-x_{i_2})^{2b}$ such that $Q_k(3,\mathbf{F}_3)$ does not contain a Ren--Xu counterexample.

From \cite{Ren_2020}, the Hilbert series for $Q_m(3,\mathbf{F}_3)$ differs from characteristic $0$ when there exists $a\in \mathbf{N}_0$ such that
\[\frac{1}{3}\leq \Bigl\{\frac{m}{3^a}\Bigr\}\leq \frac{2}{3}-\frac{1}{3^a}.\]
Notice this is equivalent to $m \pmod{3^a}$ being in $\bigl\{3^{a-1},3^{a-1}+1,\dots,2\cdot 3^{a-1}-1\bigr\}$.

\begin{Lemma} \label{noOnes}
 If $m\not\in X$, then the base $3$ representation of $m$ contains no $1$'s.
\end{Lemma}
\begin{proof}
 Suppose $m$ had the digit $1$ in the $a$-th position from the right. Then $m\pmod{3^a}$ has a~leading digit of $1$ if we choose $m\pmod{3^a}$ to be between $0$ and $3^a-1$ inclusive. However, this implies that $m \pmod{3^a}$ is in $\bigl\{3^{a-1},3^{a-1}+1,\dots,2\cdot 3^{a-1}-1\bigr\}$, so $m$ is a counterexample.
\end{proof}

\begin{Corollary}\label{evenCounter}
 If $m\not\in X$, then $m$ is even.
\end{Corollary}
\begin{proof}
 From Lemma~\ref{noOnes} $m$ has no $1$'s in its base $3$ representation, so
 \[m=\sum_{j=0}c_j3^j,\]
 where $c_j$ is $0$ or $2$. Thus $m$ must be even.
\end{proof}

\begin{Corollary}\label{noConsec}
 For all $m\not\in X$, we have $m+1\in X$.
\end{Corollary}
\begin{proof}
 By Corollary~\ref{evenCounter}, if $m\not\in X$, $m$ is even. Then $m+1$ is odd, so by the contrapositive of Corollary~\ref{evenCounter}, $m+1\in X$.
\end{proof}

Now we begin describing the degrees of Ren--Xu counterexamples.

\begin{Lemma} \label{countEx3m+3}
 If $Q_m(3,\mathbf{F}_3)_{\siv}$ has a generator in degree $3m+1$, then $m+1\in X$ and $\deg R_{m+1}=3m+3$.
\end{Lemma}
\begin{proof}
 If $m\in X$, we must have $\deg R_m<3m+1$. This implies a generator in a degree less than $3m+1$, violating Lemma~\ref{lem:3.3}. Thus $m\not\in X$, implying that $m+1\in X$ by Corollary~\ref{noConsec}.

 Because $\deg R_{m+1}<3m+4$ and $Q_{m+1}(3,\mathbf{F}_3)_\siv\subset Q_{m}(3,\mathbf{F}_3)_\siv$, we have $3m+1\leq\deg R_{m+1}< 3m+4$. By construction $3|\deg R_{m+1}$, so $\deg R_{m+1}= 3m+3$.
\end{proof}

We now introduce a few useful lemmas.

\begin{Lemma}\label{notInm+1}
 Suppose $Q_m(3,\mathbf{F}_3)_{\siv}$ has a smallest degree generator $L$ in degree $3m+1$. Assume that for all $j<m$, if $j\not\in X$, then $Q_j(3,\mathbf{F}_3)_\siv$ has a degree $3j+1$ generator. Then $Q_{m+1}(3,\mathbf{F}_3)_{\siv}$ has no nonsymmetric degree $3m+1$ or $3m+2$ element.
\end{Lemma}
\begin{proof}
 Any nonsymmetric $3m+1$ degree element in $Q_{m+1}(3,\mathbf{F}_3)_{\siv}$ must be a scalar multiple of $L$, so assume for contradiction $L$ is in $Q_{m+1}(3,\mathbf{F}_3)$. By Lemma~\ref{countEx3m+3}, $R_{m+1}=P_k^{3^a}\prod(x_{i_1}-x_{i_2})^{2b}$ is in degree $3m+3$ for natural numbers $a$, $b$, $k$. By Lemma~\ref{nonExsOnly}, $k\not\in X$ implying $P_k$ is a $3k+1$ generator of $Q_k(3,\mathbf{F}_3)_{\siv}$ using our assumption. Moreover, with any other generator in a degree less than $3m+3$ violating Lemma~\ref{lem:3.3}, $R_{m+1}$ must be generated by $L$, so $P_k^{3^a}\prod(x_{i_1}-x_{i_2})^{2b}=SL$ for some degree $2$ symmetric polynomial $S$. A degree $2$ symmetric polynomial divisible by $(x_{i_1}-x_{i_2})$ is impossible, so $S|P_k^{3^a}$ which implies either $S|P_k$ or $(x_1+x_2+x_3)|P_k$. Since $\overline{P_{k}}$ is in the $(-1)$-eigenspace of $s_{12}$, $P_{k}$ is as well and by Lemma~\ref{lem:3.1} we have $P_k=P_k'(x_1-x_2)^{2k+1}$. In both cases either $S|P_k'$ or $(x_1+x_2+x_3)|P_k'$. However, by our assumption $P_k$ is a generator, so $P_k'$ is not divisible by any nonconstant symmetric polynomial by Corollary~\ref{cor:3.2}.

 Similarly, suppose for contradiction that $K$ is a nonsymmetric element of $Q_{m+1}(3,\mathbf{F}_3)_{\siv}$ of degree $3m+2$. Since $Q_{m+1}(3,\mathbf{F}_3)_{\siv}$ has no nonsymmetric $3m+1$ degree element, $K$ must be a generator. By Lemma~\ref{lem:3.3}, $K$ is the only generator in degree less than $3m+3$, so $P_k^{3^a}\prod(x_{i_1}-x_{i_2})^{2b}$ is a symmetric polynomial multiple of $K$. However, the only symmetric polynomials of degree $1$ are multiples of $x_1+x_2+x_3$, implying $(x_1+x_2+x_3)|P_k$ which is impossible by Corollary~\ref{cor:3.2}.
\end{proof}

Note that by \cite{felder2001action}, $Q_m(3,\mathbf{Q})_{\std}$ has generators in degree $3m+1$ and $3m+2$, and by \cite{wang2022explicit}, such generators with even degree are divisible by $x_1+x_2-2x_3$. Let $\pi$ be the canonical mapping from characteristic $0$ to characteristic $3$. We then have the following lemma.

\begin{Lemma}\label{3m+2}
 Suppose $Q_m(3,\mathbf{F}_3)_{\siv}$ has a generator $L$ in degree $3m+1$. We can choose the generators of $Q_m(3,\mathbf{Q})_{\std}$ to be integer polynomials $L'$ and $(x_1+x_2-2x_3)K'$ with $\pi(K')=\pi(L')=L$. Moreover, if
 \[G=(x_1+x_2+x_3)\biggl(\frac{K'-L'}{3}\biggr)-x_3K',\]
 then
 \[\pi(G)= (x_1+x_2+x_3)\pi\biggl(\frac{K'-L'}{3}\biggr)-x_3L\]
 is a degree $3m+2$ generator for $Q_m(3,\mathbf{F}_3)_{\siv}$.
\end{Lemma}
\begin{proof}
 Let $L'$ be an arbitrary $3m+1$ degree generator of $Q_m(3,\mathbf{Q})_{\std}$ with coprime integer coefficients in the $(-1)$-eigenspace of $s_{12}$. By Lemma~\ref{lem:3.1}, $\pi(L')$ is an element of the $(-1)$-eigenspace of $s_{12}$ in $Q_m(3,\mathbf{F}_3)_{\siv}$ and if $\pi(L')$ is not a scalar multiple of $L$ then there must exist some other generator of $Q_m(3,\mathbf{F}_3)_{\siv}$ with degree less than or equal to $3m+1$. That generator and $L$ would violate Lemma~\ref{lem:3.3}, so we may set $\pi(L')=L$.

 A higher degree generator of $Q_m(3,\mathbf{Q})_{\std}$ has degree $3m+2$. With $\deg L=3m+1$ implying $m\not\in X$, $3m+2$ is even by Corollary~\ref{evenCounter}. Using \cite{wang2022explicit}, we let $(x_1+x_2-2x_3)K'$ be an arbitrary degree $3m+2$ generator for $Q_m(3,\mathbf{Q})_{\std}$ with coprime integer coefficients. Similarly, $\pi((x_1+x_2-2x_3)K')= (x_1+x_2+x_3)\pi(K')$ is an element of $Q_m(3,\mathbf{F}_3)_{\siv}$, so $\pi(K')$ is a non-symmetric polynomial of degree $3m+1$ in $Q_m(3,\mathbf{F}_3)_\siv$. Thus it must be a scalar multiple of $L$, and we may set $\pi(K')=L$.

 Let $G=(x_1+x_2+x_3)\big(\frac{K'-L'}{3}\big)-x_3K'$. Since
 \[(x_1+x_2-2x_3)K'-(x_1+x_2+x_3)L' = (x_1+x_2+x_3)\big(K'-L'\big)-3x_3K'\]
 and $\pi(K'-L')=L-L=0$, we have $G\in Q_m(3,\mathbf{Q})\cap \mathbf{Z}[x_1,x_2,x_3]$. Then
 \[\pi(G)= (x_1+x_2+x_3)\pi\biggl(\frac{K'-L'}{3}\biggr)-x_3L.\]
 If $\pi(G)$ generated by $L$, we must have $\pi(G)=c(x_1+x_2+x_3)L$ for some $c\in \mathbf{F}_3$ since $\deg\left(\pi(G)\right)=\deg(L)+1$. However, $x_1+x_2+x_3$ does not divide $x_3L$ since $L$ is a generator, so $x_1+x_2+x_3\nmid \pi(G)$. Then if $\pi(G)$ was not a generator, there must be some generator other than $L$ for $Q_m(3,\mathbf{F}_3)$ in degree less than $3m+2$ which violates Lemma~\ref{lem:3.3}. Thus, $\pi(G)$ is a generator.
\end{proof}

We aim to prove that minimum Ren--Xu counterexamples are generators and represent the only cases, where the Hilbert series of the quasi-invariants differs between characteristics $0$ and~$3$. To this end, we describe the degree of Ren--Xu counterexamples.

\begin{Example}
 We notice a ``staircase'' pattern for Ren--Xu counterexamples. The following are counterexamples for $m=3,4,5$:
 \begin{gather*}
 (x_1-x_2)^9, \qquad (x_1-x_2)^9,\qquad
 (x_1-x_2)^9\prod (x_{i_1}-x_{i_2})^2.
 \end{gather*}
 We note that since $(x_1-x_2)^9\in Q_4(3,\mathbf{F}_3)$, $(x_1-x_2)^9$ is the Ren--Xu counterexample for both~${m=3}$ and $m=4$. Moreover, the counterexample in $Q_5(3,\mathbf{F}_3)$ is the previous counterexample $(x_1-x_2)^9$ multiplied by $\prod (x_{i_1}-x_{i_2})^2$ to add the extra factor of $(x_1-x_2)^2$. In~this way the degree of counterexample stays constant for the first half of the ``staircase'' and climbs by $6$ per each increase in $m$ thereafter. Moreover, we note that $m=2,6\not\in X$, so our ``staircase'' is surrounded by non-counterexamples. One can also compute another generator for $m=3, 4, 5$ in degree $12$, $18$, and $18$ respectively. Since $9+12=6\cdot 3+3$, $9+18=6\cdot 4+3$, and $15+18=6\cdot 5 +3$, $Q_m(3,\mathbf{F}_3)_{\siv}$ is freely generated by each of these generators and $1$ by Lemma~\ref{lem:3.4}. This way we see that the upper degree generators form a complement to the lower degree ones, climbing by $6$ degrees initially and staying constant for the second half of the staircase.

 Visually, the following figure shows the degree of the generators for $Q_m(3,\mathbf{F}_3)$ with respect to $m$ were the staircase pattern and Theorem~\ref{justRens} to hold.
 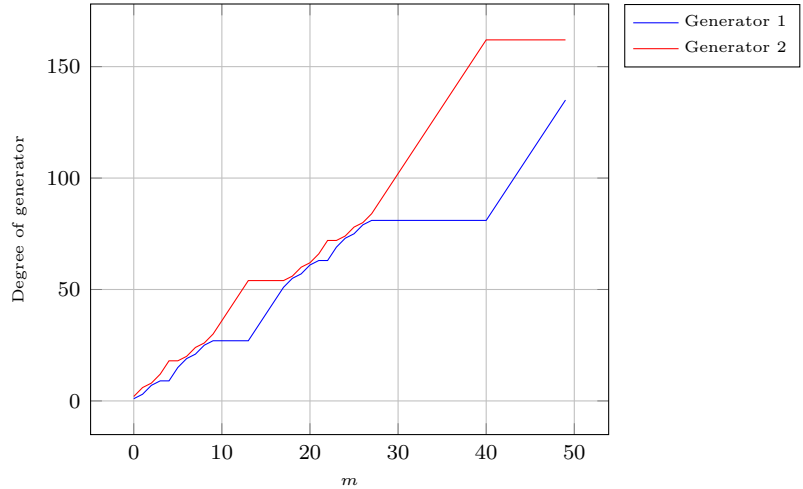
\begin{figure}[ht]
 \centering
 \begin{tikzpicture}[scale=1]
 \begin{axis}[
 tick label style={font=\scriptsize},
 xlabel={\tiny $m$},
 ylabel={\tiny Degree of generator},
 legend pos=outer north east,
 grid=major
 ]
 \addplot[color=blue] coordinates {
 (0,1)
 (1,3.0)
 (2,7)
 (3,9)
 (4,9.0)
 (5,15.0)
 (6,19)
 (7,21.0)
 (8,25)
 (9,27)
 (10,27)
 (11,27)
 (12,27)
 (13,27.0)
 (14,33.0)
 (15,39.0)
 (16,45.0)
 (17,51.0)
 (18,55)
 (19,57.0)
 (20,61)
 (21,63)
 (22,63.0)
 (23,69.0)
 (24,73)
 (25,75.0)
 (26,79)
 (27,81)
 (28,81)
 (29,81)
 (30,81)
 (31,81)
 (32,81)
 (33,81)
 (34,81)
 (35,81)
 (36,81)
 (37,81)
 (38,81)
 (39,81)
 (40,81.0)
 (41,87.0)
 (42,93.0)
 (43,99.0)
 (44,105.0)
 (45,111.0)
 (46,117.0)
 (47,123.0)
 (48,129.0)
 (49,135.0)
 };
 \addlegendentry{\tiny Generator $1$}

 \addplot[color=red] coordinates {
 (0,2)
 (1,6.0)
 (2,8)
 (3,12)
 (4,18.0)
 (5,18.0)
 (6,20)
 (7,24.0)
 (8,26)
 (9,30)
 (10,36)
 (11,42)
 (12,48)
 (13,54.0)
 (14,54.0)
 (15,54.0)
 (16,54.0)
 (17,54.0)
 (18,56)
 (19,60.0)
 (20,62)
 (21,66)
 (22,72.0)
 (23,72.0)
 (24,74)
 (25,78.0)
 (26,80)
 (27,84)
 (28,90)
 (29,96)
 (30,102)
 (31,108)
 (32,114)
 (33,120)
 (34,126)
 (35,132)
 (36,138)
 (37,144)
 (38,150)
 (39,156)
 (40,162.0)
 (41,162.0)
 (42,162.0)
 (43,162.0)
 (44,162.0)
 (45,162.0)
 (46,162.0)
 (47,162.0)
 (48,162.0)
 (49,162.0)
 };
 \addlegendentry{\tiny Generator $2$}

 \end{axis}
 \end{tikzpicture}
 \caption{Degrees of generators in characteristic $3$ with respect to $m$.} 
\end{figure}
\end{Example}

 We prove that Ren--Xu counterexamples follow this staircase pattern.

\begin{Lemma}\label{staircase}
 Let $m$ be a natural number not in $X$ and let $d$ be the largest integer such that $R_{m+1}$ lies in $Q_{m+d}(3,\mathbf{F}_3)$. Suppose that for all $k\leq m$, if $k\not\in X$, then $Q_k(3,\mathbf{F}_3)_\siv$ has a generator in degree $3k+1$. Then $R_{m+j}=R_{m+1}$ in degree $3m+3$ for $1\leq j\leq d$ and $R_{m+j}=R_{m+1}\prod (x_{i_1}-x_{i_2})^{2(j-d)}$ in degree $3m+3+6(j-d)$ for $d< j<2d$.
\end{Lemma}
\begin{proof}
 Let
 \[
 R_{m+1}=P_k^{3^a}\prod(x_{i_1}-x_{i_2})^{2b},
 \]
 where $k$ is a nonnegative integer, $a$ is a positive integer, and \smash{$b=\max\bigl\{0, \frac{2m+3-3^a(2k+1)}{2}\bigr\}$}. If $b$ is positive, the polynomial $P_k^{3^a}\prod(x_{i_1}-x_{i_2})^{2(b-1)}$ has degree less than $3m-2$ and is at least $m$-quasi-invariant since $P_k^{3^a}\prod(x_{i_1}-x_{i_2})^{2b}$ has degree less than $3m+4$. Thus $P_k^{3^a}\prod(x_{i_1}-x_{i_2})^{2(b-1)}$ is a Ren--Xu counterexample for $Q_m(3,\mathbf{F}_3)$, a~contradiction.

 In this way, we have $R_{m+1}=P_k^{3^a}$. Moreover, $Q_k(3,\mathbf{F}_3)$ must be a non-counterexample by Lemma~\ref{nonExsOnly}, so by our assumption $P_k$ is a generator. By Lemma~\ref{notInm+1}, $P_k$ is not in $Q_{k+1}(3,\mathbf{F}_3)$, so the largest power of $(x_1-x_2)$ dividing into $R_{m+1}$ must be $(x_1-x_2)^{3^a(2k+1)}$ and $m+d=\frac{3^a(2k+1)-1}{2}$ by Lemma~\ref{lem:3.1}. Then for all $1\leq j\leq d$,
 \begin{align*}
 \frac{2(m+j)+1-3^a(2k+1)}{2} \leq \frac{2(m+d)+1-3^a(2k+1)}{2}= 0.
 \end{align*}
 Thus $R_{m+j}=P_k^{3^a}=R_{m+1}$ which is indeed in degree $3m+3$ by Lemma~\ref{countEx3m+3}.

 We claim that for $d<j< 2d$, $m+j\in X$. Let $I$ be the set of integers $h$ such that a Ren--Xu counterexample for $Q_h(3,\mathbf{F}_3)$ is $P_k^{3^a}\prod(x_{i_1}-x_{i_2})^{2b}$ for some $b\in \mathbf{Z}_{\geq 0}$. By \cite{Ren_2020}, $m\in I$ if and only~if
 \[k+\frac{1}{3}\leq \frac{m}{3^a}\leq k+\frac{2}{3}-\frac{1}{3^a},\]
 which implies $I$ is $\bigl\{s,s+1,s+2,\ldots,s+3^{a-1}-1\bigr\}$ for some $s\equiv 3^{a-1} \pmod{3^a}$. Then note that $m+1\in I$, yet $m\not\in I$ since $m\not\in X$. Thus $m\equiv 3^{a-1}-1\pmod{3^a}$. Since \smash{$m+d=\frac{3^a(2k+1)-1}{2}\in I$} as well, we have $s=3^ak+3^{a-1}$, $m=3^ak+3^{a-1}-1$, and \smash{$d=\frac{3^{a-1}+1}{2}$}. Then
 \[\frac{3^a(2k+1)-1}{2}< m+j< \frac{3^{a-1}+1}{2}+ \frac{3^a(2k+1)-1}{2} = 3^ak+2\cdot 3^{a-1},\]
 so $m+j$ is in $I$ and thus in $X$.

If $R_{m+j}=P_k^{3^a}\prod(x_{i_1}-x_{i_2})^{2b}$, where \smash{$b=\frac{2(m+j)+1-3^a(2k+1)}{2}$} for $d<j< 2d$, then $m+d=\smash{\frac{3^a(2k+1)-1}{2}}$ implies $b=j-d$. Thus $R_{m+j}=P_k^{3^a}\prod(x_{i_1}-x_{i_2})^{2(j-d)}$ has degree $3m+3+6(j-d)$ as desired.
\end{proof}

In \cite{wang2022explicit}, the first author proved that generators of $Q_m(3,\mathbf{F}_p)_{\std}$ for $p>3$ lie in $\mathbf{F}_p[x_1-x_3,\allowbreak {x_2-x_3}]$ using that $\mathbf{F}_p[x_1-x_3,x_2-x_3,x_1+x_2+x_3]=\mathbf{F}_p[x_1,x_2,x_3]$. However, this is not true for $p=3$ since $x_1-x_3+x_2-x_3=x_1+x_2+x_3$ in characteristic $3$, so we instead consider the space $\mathbf{F}_3[x_1-x_3,x_2-x_3,x_3]$. From now on, we say a polynomial's degree in $x_3$ is with respect to the basis $\{x_1-x_3,x_2-x_3,x_3\}$. Moreover, in \cite{wang2022explicit} the first author defined the polynomial
\[M_d = (x_1+x_2-2x_3)^{2\{ \frac{d}{2}\}}(x_1-x_3)^{\lfloor\frac{d}{2}\rfloor}(x_2-x_3)^{\lfloor\frac{d}{2}\rfloor}\]
for natural numbers $d$ and proved that homogeneous $s_{12}$-invariant elements of $\mathbf{F}_p[x_1-x_3,x_2-x_3]/(x_1-x_2)^2$ are equal to constant multiples of $M_d$. Extending this gives that elements of $\mathbf{F}_3[x_1-x_3,x_2-x_3,x_3]/(x_1-x_2)^2$ are polynomials in $x_3$ with coefficients that are constant multiples of $M_d$. Some further nice properties of $M_d$ are the following.

\begin{Lemma}\label{M_dProps}
 For any $j,j'\in \mathbf{Z}_{\geq 0}$,
 \begin{enumerate}\itemsep=0pt
 \item[$1.$] $(x_1+x_2+x_3)M_j = M_{j+1}$ in $\mathbf{F}_3[x_1-x_3,x_2-x_3,x_3]/(x_1-x_2)^2$.
 \item[$2.$] $M_{j}M_{j'}=M_{j+j'}$ in $\mathbf{F}_3[x_1-x_3,x_2-x_3,x_3]/(x_1-x_2)^2$.
 \end{enumerate}
\end{Lemma}
\begin{proof}
1. In $\mathbf{F}_3[x_1-x_3,x_2-x_3,x_3]/(x_1-x_2)^2$, for $j\in Z_{\geq 0}$, \begin{align*}
 (x_1+x_2+x_3)M_{2j} &= (x_1+x_2+x_3)(x_1-x_3)^j(x_2-x_3)^j
 = M_{2j+1}
 \end{align*}
 and
 \begin{align*}
 (x_1+x_2+x_3)M_{2j+1} &= (x_1+x_2+x_3)^2(x_1-x_3)^j(x_2-x_3)^j \\
 &= (x_1-x_3)(x_2-x_3)M_{2j}
 = M_{2j+2}.
 \end{align*}

2. From (1), we have $M_j=(x_1+x_2+x_3)^j$ and $M_{j'}=(x_1+x_2+x_3)^{j'}$ in $\mathbf{F}_3[x_1-x_3,x_2-x_3,x_3]/(x_1-x_2)^2$, and our equality follows.
\end{proof}

This gives us intuition for the following lemmas.

\begin{Lemma}\label{monAll}
 Let $e_1$, $e_2$, and $e_3$ be the elementary symmetric polynomials for $\mathbf{F}_3[x_1,x_2,x_3]^{S_3}$ in degree $1$, $2$, and $3$ respectively. If $n$ is a natural number such that $n\not\equiv 0\pmod{3}$, for all natural numbers $j<n$ there exists a monomial $P$ in $e_1$, $e_2$, $e_3$ such that $P$ has degree $n$ and degree $j$ in~$x_3$. If $n$ is a natural number such that $n\equiv 0 \pmod{3}$, for all natural numbers $j<n-1$ there exists a monomial $P$ in $e_1$, $e_2$, $e_3$ such that $P$ has degree $n$ and degree $j$ in $x_3$.
\end{Lemma}
\begin{proof}
 We choose $e_1$, $e_2$, and $e_3$ to be
 \begin{gather*}
 e_1=x_1+x_2+x_3=(x_1-x_3)+(x_2-x_3),\\
 e_2=x_1x_2+x_1x_3+x_2x_3=(x_1-x_3)(x_2-x_3)+2((x_1-x_3)+(x_2-x_3))x_3,\\
 e_3=x_1x_2x_3=(x_1-x_3)(x_2-x_3)x_3+((x_1-x_3)+(x_2-x_3))x_3^2 + x_3^3.
 \end{gather*}

 We prove the lemma by decreasing induction on $j$.

 The base case for $n$ where $3\nmid n$ is $j=n-1$. If $j=n-1$ and $n\equiv 1 \pmod{3}$, we can let \smash{$P=e_3^{(n-1)/{3}}e_1$}. If $n\equiv 2 \pmod{3}$, we let \smash{$P=e_3^{(n-2)/{3}}e_2$}. The base case when $3|n$ is $j=n-2$,~so we can let \smash{$P=e_1e_2e_3^{(n/3)-1}$}.

 Suppose that, when $3\nmid n$, for all $j'$ such that $n>j'>j$ where $j\in \mathbf{N}$ and $0\leq j<n-1$ there exists a monomial in $e_1$, $e_2$, $e_3$ with degree $n$ and degree $j'$ in $x_3$. Suppose the same for when $3|n$ but with $n-1>j'>j$ and $j<n-2$. Then there exists a monomial $m=e_1^ae_2^be_3^c$ with degree $j+1$ in $x_3$ in $\mathbf{F}_3[x_1-x_3,x_2-x_3,x_3]/(x_1-x_2)^2$. If $b\neq 0$ we can take the monomial $e_1^{a+2}e_2^{b-1}e_3^{c}$ to be $P$ since it has degree $n$ and degree $j$ in $x_3$. If $b=0$ and $a,c>0$, then we take $P=e_1^{a-1}e_2^{b+2}e_3^{c-1}$. Finally, we are left with the cases $a,b=0$ or $b,c=0$. The former would imply \smash{$m=e_3^{\frac{n}{3}}$} is our monomial, but $3\nmid n$ would imply $m$ is not a polynomial and $3|n$ implies~$m$ has degree $j+1=n$ in $x_3$ and $j=n-1\not < n-2$. For the latter case, we have that $a=n$, so~${m=e_1^{n}}$ implies that $j+1=0$ which is below our range for $j$.
\end{proof}

\begin{Lemma}\label{symAll}
 For all $f_j\in \mathbf{F}_3$ and $n\not\equiv 0\pmod{3}$, there exists a $P\in \mathbf{F}_3[x_1,x_2,x_3]^{S_3}$ such that
 \[P= f_0M_nx_3^0+f_1M_{n-1}x_3^1+\dots + f_{n-2}M_2x_3^{n-2}+f_{n-1}M_1x_3^{n-1}\]
 in $\mathbf{F}_3[x_1,x_2,x_3]^{S_3}/(x_1-x_2)^2$. If $n\equiv 0\pmod{3}$, for all $f_j\in \mathbf{F}_3$ there exists a $P\in \mathbf{F}_3[x_1,x_2,x_3]^{S_3}$ such that
 \[P= f_0M_nx_3^0+f_1M_{n-1}x_3^1+\dots + f_{n-2}M_2x_3^{n-2}\]
 in $\mathbf{F}_3[x_1,x_2,x_3]^{S_3}/(x_1-x_2)^2$. Moreover, $P$ also satisfies the property that if it has degree $k$ in~$x_3$ in $\mathbf{F}_3[x_1,x_2,x_3]^{S_3}/(x_1-x_2)^2$, then it has degree $k$ in $x_3$ in $\mathbf{F}_3[x_1-x_3,x_2-x_3,x_3]$.
\end{Lemma}
\begin{proof}
 A weaker statement is that there exists some fixed $c_0, c_1, \dots, c_j\in\mathbf{F}_3$ such that for all $f_{j+1},f_{j+2},\dots,f_{n-1}\in \mathbf{F}_3$, there exists a symmetric polynomial
 \begin{align*}
 P&\equiv c_0M_nx_3^0+c_1M_{n-1}x_3^1+\dots + c_jM_{n-j}x_3^j\\
 &\quad +f_{j+1}M_2x_3^{n-j-1}+f_{j+2}M_1x_3^{n-j-2}+\dots+ f_{n-1}M_1x_3^{n-1}\pmod{(x_1-x_2)^2},
 \end{align*}
 when $n\not\equiv 3\pmod{3}$ and $j\in \mathbf{Z}_{\geq 0}$. A similar weaker statement can be made for the $n\equiv 0\pmod{3}$ case. We prove the statement in the lemma by induction on this $j$.

 For the base case when $n\not\equiv 0\pmod{3}$, we claim there exists coefficients $c_j\in \mathbf{F}_3$ such that the polynomial \smash{$c_0M_nx_3^0+c_1M_{n-1}x_3^1+\dots + c_{n-2}M_2x_3^{n-2}+c_{n-1}M_1x_3^{n-1}$} is in $\mathbf{F}_3[x_1,x_2,x_3]^{S_3}/(x_1-x_2)^2$. The symmetric polynomial $0$ satisfies these conditions and has degree $0$ in $x_3$. For the~base case when $n\equiv 0\pmod{3}$, we claim there exists coefficients $c_0,\dots,c_{n-2}$ such that the polynomial $c_0M_nx_3^0+c_1M_{n-1}x_3^1+\dots + c_{n-2}M_2x_3^{n-2}$ is in $\mathbf{F}_3[x_1,x_2,x_3]^{S_3}/(x_1-x_2)^2$. The~symmetric polynomial $0$ satisfies this.

 We consider the case where $n\not\equiv 0\pmod{3}$. Suppose that for all $n\geq j'>j$ there exists coefficients $c_0,\dots,c_{j'-1}$ such that for all $f_{j'},f_{j'+1},\dots, f_{n-1}$ there exists a symmetric polynomial $P$ such that
 \[P = c_0M_nx_3^0+c_1M_{n-1}x_3^1+\dots + c_{j'-1}M_{n-j'+1}x_3^{j'-1}+f_{j'}M_{n-j'}x_3^{j'}+\dots + f_{n-1}M_1x_3^{n-1}\]
 lies in $\mathbf{F}_3[x_1,x_2,x_3]^{S_3}/(x_1-x_2)^2$, where $j\in \mathbf{N}$, $0\leq j\leq n-1$. Moreover, suppose the polynomial~$P$ exists such that it has degree in $x_3$ equal to the degree in $x_3$ in $\mathbf{F}_3[x_1,x_2,x_3]^{S_3}/(x_1-x_2)^2$.

 Consider arbitrary coefficients $f_{j}, f_{j+1},\ldots, f_{n-1}$. If they are each $0$, then we can take $0$ to be our polynomial just like our base case. Otherwise, let $l$ be the greatest natural number $l\geq j$ such that $f_l\neq 0$. If $l=j$, by Lemma~\ref{monAll} there exists a monomial $m$ in $e_1$, $e_2$, $e_3$ with degree $j$ in $x_3$ and we may take $f_jm$ to be our symmetric polynomial.

 If $l>j$, by assumption there exists coefficients $c_0,c_1,\dots,c_j$ such that
 \[S=c_0M_nx_3^0+c_1M_{n-1}x_3^1+\dots + c_{j}M_{n-j}x_3^{j}+f_{j+1}M_{n-j-1}x_3^{j+1}+\dots + f_{n-1}M_1x_3^{n-1}\]
 lies in $\mathbf{F}_3[x_1,x_2,x_3]^{S_3}/(x_1-x_2)^2$. By assumption, $S$ has degree $l$ in $x_3$.

 Without loss of generality let the leading coefficient of $m$ be $M_{n-j}$, so
 \begin{gather*}
 S+(f_j-c_j)m=c_0'M_nx_3^0+c_1'M_{n-1}x_3^1+\dots + c_{j-1}'M_{n-j+1}x_3^{j-1}\\
 \hphantom{S+(f_j-c_j)m=}{} +f_{j}M_{n-j}x_3^{j}+\dots + f_{n-1}M_1x_3^{n-1}\
 \end{gather*}
 for some coefficients $c_0',c_1',\ldots, c_{j-1}'$. Moreover, $S+(f_j-c_j)m$ is still a symmetric polynomial and $m$ has degree $j$ in $x_3$ while $S$ has degree $l$, so $S+(f_j-c_j)m$ has degree $l$ as desired.

 An identical argument holds for $n\equiv 0\pmod{3}$.
\end{proof}

Now we have the tools to prove $m\not\in X$ implies $m+1$ begins our staircase.

\begin{Lemma}\label{3m+6}
 Suppose that for all $k\leq m$, if $k\not\in X$ then $Q_k(3,\mathbf{F}_3)$ has a $3k+1$ degree generator, where $m$ is a natural number. Then if $Q_m(3,\mathbf{F}_3)_{\siv}$ has a generator in degree~${3m+1}$, $Q_{m+1}(3,\mathbf{F}_3)_{\siv}$ has a generator in degree $3m+6$.
\end{Lemma}
\begin{proof}
 By Lemma~\ref{3m+2}, the generators for $Q_m(3,\mathbf{F}_3)_{\siv}$ are
 \[\biggl((x_1+x_2+x_3)\pi\biggl(\frac{A'-B'}{3}\biggr)-x_3B\biggr)(x_1-x_2)^{2m+1}\]
 in degree $3m+2$, and
 \[B(x_1-x_2)^{2m+1}\]
 in degree $3m+1$, where $(x_1-x_2)^{2m+1}(x_1+x_2-2x_3)A'$ and $(x_1-x_2)^{2m+1}B'$ are the generators of $Q_m(3,\mathbf{Q})_{\std}$, $B$ is an $s_{12}$-invariant polynomial, and $\pi(A')=\pi(B')=B$.

 For the greater degree generator, let \smash{$C=\big((x_1+x_2+x_3)\pi\big(\frac{A'-B'}{3}\big)-x_3B\big)$}. We would like to show there exists symmetric polynomials $P$ and $Q$ in degree $4$ and $5$ respectively such that
 \[PC+QB \equiv 0\pmod{(x_1-x_2)^2}.\]
 Since \smash{$\frac{PC+BQ}{(x_1-x_2)^2}$} is still $s_{12}$-invariant, this would then imply $(PC+QB)(x_1-x_2)^{2m+1}\in Q_{m+1}(3,\mathbf{F}_3)$ by Lemma~\ref{lem:3.1}. Consider writing
 \[P=f_0M_4x_3^0+f_1M_{3}x_3^1 + f_{2}M_2x_3^{2}+f_{3}M_1x_3^{3}\]
 and
 \[Q=h_0M_5x_3^0+h_1M_4x_3^1+h_2M_3x_3^2+h_3M_2x_3^3+h_4M_1x_3^4\]
 for arbitrary $f_j$ and $h_j$. By Lemma~\ref{symAll}, we know that for any choice of $f_j$ and $h_j$, we have $P,Q\in \mathbf{F}_3[x_1,x_2,x_3]^{S_3}/(x_1-x_2)^2$.

 We claim that \smash{$B|\pi\big(\frac{A'-B'}{3}\big)$} in $\mathbf{F}_3[x_1,x_2,x_3]/(x_1-x_2)^2$. By \cite{wang2022explicit}, $A'$ and $B'$ are both polynomials in the variables $(x_1-x_2)^2$ and $(x_1-x_3)(x_2-x_3)$. Moreover, by Lemma~\ref{notInm+1}, $(x_1-x_2)^2\nmid B$ so $B\equiv cM_{m}\pmod{(x_1-x_2)^2}$ for some $c\in \mathbf{F}_3$ such that $c\neq 0 $. Similarly, we know $\smash{\pi\big(\frac{A'-B'}{3}\big)} \equiv c'M_{m}\pmod{(x_1-x_2)^2}$ for some $c'\in \mathbf{F}_3$. Thus we have $\smash{\pi\big(\frac{A'-B'}{3}\big)} = dB$, where $d=\frac{c'}{c}$.

 We use Lemma~\ref{M_dProps} to expand $PC+BQ$ in $\mathbf{F}_3[x_1,x_2,x_3]/(x_1-x_2)^2$,
 \begin{align*}
 PC+QB&=\biggl(h_0M_5B+f_0(x_1+x_2+x_3)M_4\pi\biggl(\frac{A'-B'}{3}\biggr)x_3^0\biggr) \\
 &\quad + \sum _{j=1}^{3}\biggl(h_{j}M_{5-j}Bx_3^j+f_{j}(x_1+x_2+x_3)M_{4-j}\pi\biggl(\frac{A'-B'}{3}\biggr)x_3^{j}\\
 &\quad\hphantom{+ \sum _{j=1}^{3}\biggl(}{}- f_{j-1}M_{5-j}Bx_3^j\bigg)+ h_4M_1Bx_3^4-f_3M_1Bx_3^4 \\
 &=\biggl(h_0B+f_0\pi\biggl(\frac{A'-B'}{3}\biggr)\biggr)M_5 \\
 &\quad + \sum _{j=1}^{3}\biggl(\biggl((h_{j} - f_{j-1})B+f_{j}\pi\biggl(\frac{A'-B'}{3}\biggr)\biggr)M_{5-j}x_3^j\biggr) + (h_4-f_3)M_1Bx_3^4. \\
 &=(h_0+f_0d)BM_5 + \sum _{j=1}^{3}\biggl((h_{j} - f_{j-1})+f_{j}d\biggr)BM_{5-j}x_3^j + (h_4-f_3)M_1Bx_3^4.
 \end{align*}
 Letting $h_j$ be arbitrary for $j>0$, set $f_3=h_4$, $f_{j-1} = h_{j}+f_jd$ for $0<j<3$ and set $h_0=-f_0d$. This makes the expression $PC+QB=0$.

 We claim $Q_{m+1}(3,\mathbf{F}_3)_{\siv}$ has a degree $3m+3$ generator, namely $R_{m+1}$. From Lemma~\ref{notInm+1}, $Q_{m+1}(3,\mathbf{F}_3)_{\siv}$ has no degree $3m+1$ or $3m+2$ generator, so it has no generators in degree less than $3m+3$. By Lemma~\ref{countEx3m+3}, $R_{m+1}$ is in degree $3m+3$ so it must be a generator. Without loss of generality, we let
 \[R_{m+1}=\left((x_1+x_2+x_3)C + SB\right)(x_1-x_2)^{2m+1},\]
 where $S$ is a degree $2$ symmetric polynomial.

 If $(PC+QB)(x_1-x_2)^{2m+1}$ were generated by $R_{m+1}$, there would exist a symmetric polynomial~$I$ such that $IR_{m+1}=(PC+QB)(x_1-x_2)^{2m+1}$. This implies $(I(x_1+x_2+x_3)-P)C+(IS-Q)B=0$. If $I(x_1+x_2+x_3)-P\neq 0$ or $IS-Q\neq 0$, there is a relation on $C$ and $B$ over $\mathbf{F}_3[x_1,x_2,x_3]^{S_3}$, but $C(x_1-x_2)^{2m+1}$ and $B(x_1-x_2)^{2m+1}$ are generators of $Q_m(3,\mathbf{F}_3)$. Thus we must have $P= I(x_1+x_2+x_3)$, so $(x_1+x_2+x_3)|P$. Now we consider the symmetric polynomials $P'=P+e_2^2+e_2e_1^2+e_1^4$ and $Q' = Q+e_3e_1^2+(-d-1)e_2^2e_1-de_1^3e_2+(-d+1)e_1^5$. In $F_3[x_1-x_3,x_2-x_3,x_3]/(x_1-x_2)^2$, we get that
 \[P'=f_0M_4x_3^0+f_1M_{3}x_3^1 + (f_{2}+1)M_2x_3^{2}+f_{3}M_1x_3^{3}\]
 and
 \[Q'=h_0M_5x_3^0+h_1M_4x_3^1+(h_2-d)M_3x_3^2+(h_3+1)M_2x_3^3+h_4M_1x_3^4.\]
 Then $f_2+1 = (h_3+f_3d)+1=(h_3+1)+f_3d$, $f_1=h_2+f_2d = (h_2-d)+(f_2+1)d$, and the rest of the equations necessary for $P'C+Q'B\equiv 0\pmod{(x_1-x_2)^2}$ are the same as $PC+QB\equiv 0\pmod{(x_1-x_2)^2}$. Thus $P'C+Q'B\equiv 0\pmod{(x_1-x_2)^2}$. Moreover, $(x_1+x_2+x_3)$ divides into $P+e_2e_1^2+e_1^4$ but not $e_2^2$, so $(x_1+x_2+x_3)\nmid P'$. We have shown that if ${(PC+QB)(x_1-x_2)^{2m+1}}$ is generated by $R_{m+1}$, then $(x_1+x_2+x_3)|P$, implying ${(P'C+Q'B)(x_1-x_2)^{2m+1}}$ is not generated by $R_{m+1}$. If $(P'C+Q'B)(x_1-x_2)^{2m+1}$ is not a generator, then whatever generates it violates~Lemma~\ref{lem:3.3}, so $(P'C+Q'B)(x_1-x_2)^{2m+1}$ is indeed a degree $3m+6$ generator of~$Q_{m+1}(3,\mathbf{F}_3)$.
 \end{proof}

Now we prove that if $R_{m+1}$ begins our staircase, then it is the lower degree generator for the first half of the staircase.

\begin{Lemma}\label{higerGen}
 Let $m\not\in X$ for some natural number $m$. Suppose $R_{m+1}$ is a degree $3m+3$ generator of $Q_{m+1}(3,\mathbf{F}_3)_{\siv}$ and $L$ is another generator in degree $3m+6$. Further, let $R_{m+1}$ lie in $Q_{m+d}(3,\mathbf{F}_3)$, where $d$ is maximal. Then $L\prod(x_{i_1}-x_{i_2})^{2(j-1)}$, $R_{m+1}$, and $1$ freely generate $Q_{m+j}(3,\mathbf{F}_3)_{\siv}$ for $1\leq j \leq d$.
\end{Lemma}
\begin{proof}
 As a generator, $L$ lies in a copy of $\siv$ and is divisible by $(x_1-x_2)^{2(m+1)+1}$ by Lemma~\ref{lem:3.1}. Since $L\prod(x_{i_1}-x_{i_2})^{2(j-1)}$ is divisible by $(x_1-x_2)^{2(m+j)+1}$, by the second part of Lemma~\ref{lem:3.1}, \smash{$L\prod(x_{i_1}-x_{i_2})^{2(j-1)}$} is in \smash{$Q_{m+j}(3,\mathbf{F}_3)_{\siv}$}. If $L\prod(x_{i_1}-x_{i_2})^{2(j-1)}$ is not a~generator, $R_{m+1}$ must generate \smash{$L\prod(x_{i_1}-x_{i_2})^{2(j-1)}$}, implying a relation between $R_{m+1}$ and~$L$. Thus $L\prod(x_{i_1}-x_{i_2})^{2(j-1)}$ is indeed a generator.

 Moreover, $3m+3+3m+6j=6(m+j)+3$ so by Lemma~\ref{lem:3.4}, $L\prod(x_{i_1}-x_{i_2})^{2(j-1)}$ and $R_{m+1}$ generate $Q_{m+j}(3,\mathbf{F}_3)_{\siv}$.
\end{proof}

Next, we prove that, for all consecutive spaces of quasi-invariants in the second half of the staircase, the lower degree generator is $\prod (x_{i_1}-x_{i_2})^{2}$ times the previous lower degree generator.

\begin{Lemma}\label{lowerGen}
 Let $m\not\in X$ for some natural number $m$. Suppose $R_{m+1}$ is a degree $3m+3$ generator of $Q_m(3,\mathbf{F}_3)$ and $L$ is another generator in degree $3m+6$. Let $R_{m+1}$ lie in $Q_{m+d}(3,\mathbf{F}_3)$, where $d$ is maximal. Further, let $L$ have degree at most $5$ in $x_3$. Then for all $d\leq j< 2d$, $Q_{m+j}(3,\mathbf{F}_3)_{\siv}$ is freely generated by a generator in degree $3m+6d$, $R_{m+1}\prod(x_{i_1}-x_{i_2})^{2(j-d)}$ in degree $3m+6(j-d)+3$, and $1$.
\end{Lemma}

\begin{proof} We proceed with induction.

 The generator $R_{m+1}$ of $Q_{m+d}(3,\mathbf{F}_3)$ is in degree $3m+3=3m+6(d-d)+3$, and from Lemma~\ref{higerGen} a second generator is $L\prod(x_{i_1}-x_{i_2})^{2(d-1)}$ in degree $3m+6d$. Moreover, these are the only generators so the claim is true for $j=d$.

 Let $k$ be a natural number with $d<k<2d$ and suppose $Q_{m+j}(3,\mathbf{F}_3)_{\siv}$ has a generator in degree $3m+6d$ and degree $3m+6(j-d)+3$ for all $d\leq j < k$, where this upper degree generator is a polynomial of degree at most $5$ in $x_3$ and is not generated by $R_{m+1}$. Consider $Q_{m+k}(3,\mathbf{F}_3)_{\siv}$. We know $R_{m+1}\prod(x_{i_1}-x_{i_2})^{2(k-d-1)}$ is an element of $Q_{m+k-1}(3,\mathbf{F}_3)_{\siv}$ of degree $3m+6(k-d-1)+3$ by Lemma~\ref{lem:3.1}. Since $k-1<k$, our inductive hypothesis implies $R_{m+1}\prod(x_{i_1}-x_{i_2})^{2(k-d-1)}$ is a generator for $Q_{m+k-1}(3,\mathbf{F}_3)_{\siv}$.

 Let $T$ be the degree $3m+6d$ generator for $Q_{m+k-1}(3,\mathbf{F}_3)_{\siv}$ with degree $5$ in $x_3$. We write $R_{m+1}\prod(x_{i_1}-x_{i_2})^{2(k-d-1)} = R_{m+1}'(x_1-x_2)^{2(m+k-1)+1}$ and $T=T'(x_1-x_2)^{2(m+k-1)+1}$ for $s_{12}$ invariant polynomials $R_{m+1}'$ and $T'$. If $o = m+4k-6d-2$ and $r = m+6d-2k+1$, then $\deg R_{m+1}' =o$ and $\deg T' = r$. We want to find a degree $r-o$ symmetric polynomial $P$ such that
 \[-PR_{m+1}'+T'\equiv 0 \pmod{(x_1-x_2)^2}.\]

 We claim that $R_{m+1}'$ has degree $0$ in $x_3$. This is because $R_{m+1}\prod (x_{i_1}-x_{i_2})^{2(k-d-1)}=P_l^{3^a}\prod (x_{i_1}-x_{i_2})^{2(k-d-1)}$ as we proved in Lemma~\ref{staircase}. Since $P_l$ is the map of the generator of $Q_l(3,\mathbf{Q})$ into characteristic $3$, $P_l$ must be constant in the variable $x_3$. We can see ${\prod (x_{i_1}-x_{i_2})^{2(k-d-1)}}$ is also constant in $x_3$, so $R_{m+1}$ and $R_{m+1}'$ are constant in $x_3$.

 Having assumed that $T'$ is at most degree $5$ in $x_3$,
 \[T' = t_0M_rx_3^0+t_1M_{r-1}x_3^1+t_2M_{r-2}x_3^2+t_3M_{r-3}x_3^3+t_4M_{r-4}x_3^4+t_5M_{r-5}x_3^5\]
 and
 \[R_{m+1}'=aM_{o}\]
 for coefficients $t_j$ and $a$ in $\mathbf{F}_3$. Since $R_{m+1}$ is not in $Q_{m+d+1}(3,\mathbf{F}_3)_{\siv}$, we have $a\neq 0$. We~let
 \begin{gather*}
 P = \frac{t_0}{a}M_{r-o}x_3^0+\frac{t_1}{a}M_{r-o-1}x_3^1+\frac{t_2}{a}M_{r-o-2}x_3^2+\frac{t_3}{a}M_{r-o-3}x_3^3\\
 \hphantom{P =}{}
 +\frac{t_4}{a}M_{r-o-4}x_3^4+\frac{t_5}{a}M_{r-o-5}x_3^5,
 \end{gather*}
 so that $T'-PR_{m+1}'\equiv 0\pmod{(x_1-x_2)^2}$ by Lemma~\ref{M_dProps}. Since $\deg(P) = r-o = 12d-6k+3\geq 9>7$, by Lemma~\ref{symAll} such a symmetric polynomial $P$ is attainable with $P$ having degree at most degree $5$ in $x_3$.
 Since $T'$ also has at most degree $5$ in $x_3$ and $R_{m+1}'$ has degree $0$, $(-PR_{m+1}'+T')$ has at most degree $5$ in $x_3$. Letting $U=\bigl(-PR_{m+1}'+T'\bigr)(x_1-x_2)^{2(m+k-1)+1}$, we have $U$ is in $Q_{m+k}(3,\mathbf{F}_3)$ with degree $3m+6d$ and since $\bigl(-PR_{m+1}'\bigr)(x_1-x_2)^{2(m+k-1)+1}$ is generated by~$R_{m+1}$ and $T$ is not, $U$ is not generated by $R_{m+1}$. Finally, we also have $R_{m+1}\prod (x_{i_1}-x_{i_2})^{2(k-d)}$ is in $Q_{m+k}(3,\mathbf{F}_3)_{\siv}$ with degree $3m+6(k-d)+3$. Thus what is left is to prove is $R_{m+1}\prod (x_{i_1}-x_{i_2})^{2(k-d)}$ and $\bigl(-PR_{m+1}'+T'\bigr)(x_1-x_2)^{2(m+k-1)+1}$ are generators for $Q_{m+k}(3,\mathbf{F}_3)$.

 Assume for sake of contradiction that $U$ and $R_{m+1}\prod (x_{i_1}-x_{i_2})^{2(k-d)}$ are not both generators. If $R_{m+1}\prod (x_{i_1}-x_{i_2})^{2(k-d)}$ is a generator, then any other generator must be of at least degree $3m+6d$ by Lemma~\ref{lem:3.3}. Yet $U$ is not generated by $R_{m+1}\prod (x_{i_1}-x_{i_2})^{2(k-d)}$ since it is not generated by $R_{m+1}$. Thus $U$ must be a generator.

 Next, we consider if $R_{m+1}\prod (x_{i_1}-x_{i_2})^{2(k-d)}$ is not a generator. For $R_{m+1}\prod (x_{i_1}-x_{i_2})^{2(k-d)}$ to not be a generator there must be a generator in a degree less than $3m+6(k-d)+3$. Let it be $G$, and by Lemma~\ref{lem:3.3}, any other generator must have degree greater than $3m+6d$. Thus $U$ is not a generator, so $U$ and $R_{m+1}\prod (x_{i_1}-x_{i_2})^{2(k-d)}$ are both generated by $G$ and specifically $U=QG$ and $R_{m+1}\prod (x_{i_1}-x_{i_2})^{2(k-d)}=SG$ for symmetric polynomials $P$ and $Q$. Moreover, $R_{m+1}\prod (x_{i_1}-x_{i_2})^{2(k-d-1)}$ is the lowest degree generator for $Q_{m+k-1}(3,\mathbf{F}_3)_{\siv}$, so $G=CR_{m+1}\prod (x_{i_1}-x_{i_2})^{2(k-d-1)}$ for a symmetric polynomial $C$. This implies $C|\prod(x_{i_1}-x_{i_2})^2$, and $G$ is not a scalar multiple of $R_{m+1}\prod (x_{i_1}-x_{i_2})^{2(k-d)}$, so $C$ is a constant. We then have~$U$ is a constant multiple of $QR_{m+1}\prod (x_{i_1}-x_{i_2})^{2(k-d-1)}$, so $U$ is generated by $R_{m+1}$ which is a~contradiction.

 Thus $U$ and $R_{m+1}\prod (x_{i_1}-x_{i_2})^{2(k-d)}$ are each generators and together with $1$ they freely generate $Q_m(3,\mathbf{F}_3)_{\siv}$ by Lemma~\ref{lem:3.4}.
\end{proof}

Finally, we show that after the staircase completes, the next space of quasi-invariants has no counterexamples.

\begin{Lemma}\label{backNon}
 Let $Q_{m-1}(3,\mathbf{F}_3)_{\siv}$ have generators $K$ in degree $3m-3$ and $T$ in degree $3m$ such that $K$ is not in $Q_{m}(3,\mathbf{F}_3)_{\siv}$. If $m$ is even, then $Q_{m}(3,\mathbf{F}_3)_{\siv}$ is freely generated by a generator in degree $3m+1$, $3m+2$, and $1$.
\end{Lemma}
\begin{proof}
 Suppose for the sake of contradiction that $Q_{m}(3,\mathbf{F}_3)_{\siv}$ has a generator $U$ in degree $3m-1$ or $3m-2$. Then since $U$ is also in the $-1$ $s_{12}$ eigenspace of $Q_{m-1}(3,\mathbf{F}_3)_{\siv}$, $U$ must be generated by $K$ over $\mathbf{F}[x_1,x_2,x_3]^{S_3}$. Yet $K$ being divisible by a symmetric polynomial violates Corollary~\ref{cor:3.2}.

 Suppose for the sake of contradiction that $Q_{m}(3,\mathbf{F}_3)_{\siv}$ has a generator in degree $3m$. Without loss of generality let that generator be $T$. From \cite{felder2001action}, we can let $L'$ be a degree $3m+1$ generator of $Q_m(3,\mathbf{Q})_{\std}$ with coprime integer coefficients. Then $\pi(L')\in Q_m(3,\mathbf{F}_3)_{\siv}$, so~$\pi(L')$ must be generated by $T$ since any other generator in degree less than degree $3m+1$ would violate Lemma~\ref{lem:3.3}. Moreover, the only degree $1$ symmetric polynomials are constant multiples of $x_1+x_2+x_3$, so we can assume without loss of generality that
 \[\pi(L')=(x_1+x_2+x_3)T.\]
 Note that from \cite{wang2022explicit} all generators of $Q_m(3,\mathbf{Q})_{\std}$ must lie in $\mathbf{Q}[x_1-x_3,x_2-x_3]$. Thus $(x_1+x_2+x_3)T\in\mathbf{F}_3[x_1-x_3,x_2-x_3]$ and so $T\in \mathbf{F}_3[x_1-x_3,x_2-x_3]$.

 We also have $T=(x_1-x_2)^{2m+1}T'$ for some $s_{12}$-invariant polynomial $T'$. Thus by the fundamental theorem of symmetric polynomials $T'\in \mathbf{F}_3[(x_1-x_3)(x_2-x_3),x_1+x_2+x_3]$. Note that $\deg T' = 3m-2m-1=m-1$ and $m$ is even, so $T'$ has an odd degree. However, since it is generated by $(x_1-x_3)(x_2-x_3)$ and $x_1+x_2+x_3$, we must have $(x_1+x_2+x_3)|T'$. This gives a contradiction because $T$ is a generator.
\end{proof}

Finally, we have the lemmas to prove Theorem~\ref{justRens}.

\begin{proof}[Proof of Theorem~\ref{justRens}]
 We prove this using induction on $m$.

 The generators for $Q_0(3,\mathbf{F}_3)_{\siv}$ are $x_1-x_2$ and $x_3(x_1-x_2)$. These generators are in degree $3\cdot 0 +1$ and $3\cdot 0 +2$ so the theorem is true for the base case.

 Assume the claim is true when $m<j$ for some $j\in \mathbf{N}$. Consider the space $Q_j(3,\mathbf{F}_3)_{\siv}$. Let $t$ be the largest natural number less than $j$ such that $t\not\in X$. By the inductive hypothesis, $Q_t(3,\mathbf{F}_3)$~has a generator in degree $3t+1$ and $3t+2$. By Lemma~\ref{3m+2}, we may let the generators~be
 \[\biggl(x_1+x_2+x_3)\pi\biggl(\frac{A'-B'}{3}\biggr)-x_3B\biggr)(x_1-x_2)^{2t+1}\qquad
 \text{and}\qquad
 B(x_1-x_2)^{2t+1},\]
 where $(x_1-x_2)^{2t+1}(x_1+x_2-2x_3)A'$ and $(x_1-x_2)^{2t+1}B'$ are generators for $Q_t(3,\mathbf{Q})_{\std}$ and $\pi(A')=\pi(B')=B$. From Lemma~\ref{3m+6}, $Q_{t+1}(3,\mathbf{F}_3)_{\siv}$ is generated by a generator in degree ${3t+6}$ and $3t+3$. Moreover, $R_{t+1}$ is the $3t+3$ degree generator by Lemma~\ref{countEx3m+3}. Let $L$ be the degree $3t+6$ generator. Suppose $R_{t+1}$ lies in $Q_{t+d}(3,\mathbf{F}_3)$, but not $Q_{t+d+1}(3,\mathbf{F}_3)$, where $d$ is a~natural number.

 First, we consider when $t+d\geq j\geq t+1$. By Lemma~\ref{higerGen}, $Q_j(3,\mathbf{F}_3)_{\siv}$ has generators $R_{t+1}$ and $L\prod (x_{i_1}-x_{i_2})^{2(j-t-1)}$. Note that $R_{t+1}=R_{j}$ by Lemma~\ref{staircase}, and further $\deg\big(L\prod (x_{i_1}-x_{i_2})^{2(j-t-1)}\big) + \deg(R_{t+1}) = (6(j-t-1) + 3t+6)+3t+3=6j+3$. By Lemma~\ref{lem:3.4}, we then have that $R_{t+1}$ and \smash{$L\prod (x_{i_1}-x_{i_2})^{2(j-t-1)}$} generate $Q_j(3,\mathbf{F}_3)_{\siv}$.

 Next, we consider the case where $t+2d-1\geq j\geq t+d+1$. Notice that by our construction in Lemma~\ref{3m+6}, we can choose $L$ such that it has at most degree $5$ in $x_3$. Thus we can apply Lemma~\ref{lowerGen}, which gives us that $Q_j(3,\mathbf{F}_3)_{\siv}$ is generated by $R_{t+1}\prod (x_{i_1}-x_{i_2})^{2(j-t-d)}$ and a generator in degree $3t+6d$. Note that $R_{t+1}\prod (x_{i_1}-x_{i_2})^{2(j-t-d)}$ is a constant multiple of $R_j$ by Lemma~\ref{staircase}. Moreover, the sum of their degrees is $3t+6(j-t-d)+3+3t+6d=6j+3$ as desired.

 Finally, we consider if $j=t+2d$. Note that by Lemma~\ref{lowerGen}, $Q_{t+2d-1}(3,\mathbf{F}_3)$ has a generator in degree $3t+6d$ and $3t+6(d-1)+3$. The degree $3t+6(d-1)+3$ generator is $R_{t+1}\prod (x_{i_1}-x_{i_2})^{2(d-1)}$, and $R_{t+1}$ is divisible by $(x_1-x_2)^{2(t+d)+1}$, where $d$ is maximal, so $R_{t+1}\prod (x_{i_1}-x_{i_2})^{2(d-1)}$ does not lie in $Q_{t+2d}(3,\mathbf{F}_3)$. Moreover, $Q_{t}(3,\mathbf{F}_3)$ is a non Ren--Xu counterexample, so $t$ must be even by Lemma~\ref{evenCounter}. Then $t+2d$ is even as well, so by Lemma~\ref{backNon}, $Q_{t+2d}(3,\mathbf{F}_3)$ has a generator in degree $3(t+2d)+1$ and $3(t+2d)+2$.

 Now we claim we have exhausted all cases. If we had $j>t+2d$, since we just showed ${t+2d\not\in X}$, we would not have chosen $t$ to be the largest natural number less than $j$ not in~$X$.
\end{proof}

\begin{Remark}
 We can compute the degrees of generators of $Q_m(3,\mathbf{F}_3)_{\siv}$ explicitly. If~$m$ has no digits $1$ in its base $3$ representation, then the generators have degree $3m+1$ and $3m+2$. Otherwise, the lower degree generator is~$R_m$. We can deduce the minimal degree of the \mbox{Ren--Xu} counterexamples in $Q_m(3,\mathbf{F}_3)$: Let $a$ be the greatest natural number such that the $a$-th term from the right in the base $3$ representation of~$m$ is~$1$. Then if \smash{$\bigl\lceil\frac{\lceil \frac{m}{3^a}\rceil-1}{2}\bigr\rceil=k$}, \mbox{a~minimal} deg\-ree~Ren--Xu counterexample is $P_k^{3^a}\prod(x_{i_1}-x_{i_2})^{2b}$, where
 \[
 b=\max\left\{\frac{2m+1-3^a(2k+1)}{2},0\right\}.
 \]
 The~deg\-rees~of the generators are then $3^a(2k+1)+6b$ and $6m+3-3^a(2k+1)-6b$.
\end{Remark}

\section[Representations of S\_3 in Q\_m(3,F)]{Representations of $\boldsymbol{S_3}$ in $\boldsymbol{Q_m(3,\mathbf{F}_3)}$}\label{sec:reps}

Now that we have a complete picture of $Q_m(3,\mathbf{F}_3)_\siv$, we consider generators that generate the other indecomposable modules of $S_3$. We start with $\trix$, which behaves very similarly to $\siv$.

 \begin{Proposition}\label{prop:trix}
 Suppose that for all $j\leq m$, $Q_j(3,\mathbf{F}_3)_\siv$ has generators in degree $d$ and $6j+3-d$ respectively for some $d$. If $K$, $L$ are distinct generators of $Q_m(3,\mathbf{F}_3)_\siv$ then there are two other homogeneous generators $K_1$, $L_1$ of $Q_m(3,\mathbf{F}_3)$ in the same degrees as~$K$,~$L$, respectively such that as a representation of $S_3$, $K_1$ generates a copy of $\trix$ containing $K$ and $L_1$ generates a copy of $\trix$ containing $L$. Moreover, there are no relations between $K_1$, $L_1$ over the symmetric polynomials, and there are no other generators of $Q_m(3,\mathbf{F}_3)$ that generate a copy of $\trix$.
 \end{Proposition}

 \begin{proof}
 We prove this by induction on $m$. For the base case $m=0$, note that by Example~\ref{deg1polys}, for ${K=x_1-x_2}$ we have that $K_1=x_1$ satisfies the desired conditions. Similarly, for $L=(x_1-x_2)x_3$, we have that $L_1=x_1(x_2+x_3)$ satisfies the desired conditions. These two are independent over the symmetric polynomials, as a relation between them would imply a relation between 1 and~${x_2+x_3}$.\looseness=-1

 For the inductive step, let $K'$, $L'$ be the generators of $Q_{m-1}(3,\mathbf{F}_3)_\siv$ and let $K'_1$, $L'_1$ be the corresponding generators of $Q_{m-1}(3,\mathbf{F}_3)$. Without loss of generality, we can choose~$K'_1$,~$L'_1$ to be $s_{23}$-invariant with $(1-s_{12})K'_1=K'$, $(1-s_{12})L'_1=L'$ (similar to in the base case). Let~$K$,~$L$ be generators of $Q_m(3,\mathbf{F}_3)_\siv$. Then since $K,L\in Q_{m-1}(3,\mathbf{F}_3)_\siv$, we can write $K=P_1K'+Q_1L'$, $L=P_2K'+Q_2L'$ for symmetric polynomials $P_1$, $P_2$, $Q_1$, $Q_2$. Then it follows that $K_1:=P_1K'_1+Q_1L'_1$, $L_1:=P_2K'_1+Q_2L'_1$ each generate a copy of $\trix$ that contains $K$, $L$, respectively. Moreover, if there is some relation $P_3K_1+Q_3L_1=0$ for symmetric polynomials $P_3$, $Q_3$, then applying $1-s_{12}$ to this equation would yield $P_3K+Q_3L=0$, which violates Lemma~\ref{lem:3.4}.

 Next, we show that $K_1$, $L_1$ are $m$-quasi-invariants. As the computations are the same for both polynomials, we give the proof only for $K_1$. First, note that $(1-s_{23})K_1=0$ since both $K'_1$, $L'_1$ are $s_{23}$-invariant. Next, note that $(1-s_{12})K_1=K$ is divisible by $(x_1-x_2)^{2m+1}$ by Lemma~\ref{lem:3.1}. Finally, note that since $K_1$ is $s_{23}$-invariant, we have
 \[(1-s_{13})K_1=s_{23}(s_{23}-s_{23}s_{13})K_1=s_{23}(1-s_{23}s_{13}s_{23})K_1=s_{23}(1-s_{12})K_1\]
 is divisible by $s_{23}(x_1-x_2)^{2m+1}=(x_1-x_3)^{2m+1}$.

 Note that $K_1$, $L_1$ are the minimal degree polynomials such that $(1-s_{12})K_1$, $(1-s_{12})L_1$ are symmetric polynomial multiples of $K$, $L$, respectively, so they cannot be generated by any other generators and thus must be generators themselves. Then assume for contradiction that there is some other generator $T$ of $Q_m(3,\mathbf{F}_3)$ that generates a copy of $\trix$. Then $(1-s_{12})T$ is contained in a copy of $\siv$ and is $s_{12}$-antiinvariant, so we can write $(1-s_{12})T=S_1K+S_2L$ for symmetric polynomials $S_1$, $S_2$. Then $T$,
 $S_1K_1+S_2L_1$ generate copies of $\trix$ with the same $\siv$ submodule, so they generate a copy of
 \[(\trix\oplus\trix)/\siv\cong\trix\oplus\triv. \]
 Thus $T$ is generated by $K_1$, $L_1$, $1$, and is not a generator itself.
 \end{proof}

 \begin{Corollary}\label{cor:5gens}
 The generators $1$, $K$, $K_1$, $L$, $L_1$ of $Q_m(3,\mathbf{F}_3)$ defined in Proposition~$\ref{prop:trivsign}$, Theorem~$\ref{justRens}$ and Proposition~$\ref{prop:trix}$ have no relations between them over the symmetric polynomials.
 \end{Corollary}

 \begin{proof}
 Let
 \[P_1+P_2K+P_3L+P_4K_1+P_5L_1=0\]
 for symmetric polynomials $P_1,\dots,P_5$. Then apply $1+s_{12}$ to the equation to yield
 \[2P_1+P_4(2K_1-K)+P_5(2L_1-L)=0\]
 since $K$, $L$ are $s_{12}$-antiinvariant. Next, apply $1-s_{23}$ to this equation to yield
 \[P_4(s_{23}-1)K+P_5(s_{23}-1)L=0.\]
 Note that $(s_{23}-1)K$ generates the same copy of $\siv$ as $K$, since $s_{23}-1$ acts bijectively on $\sign$ (and similarly for $L$). So a relation between $(s_{23}-1)K$, $(s_{23}-1)L$ is equivalent to a~relation between $K$, $L$, which cannot exist by Lemma~\ref{lem:3.4}. So we have $P_4=P_5=0$.

 Now, the result follows from Lemma~\ref{lem:3.4}.
 \end{proof}

 \begin{Remark}\label{rem:sign}
 In the non-modular case, one has that the polynomial $\prod_{i_1<i_2}(x_{i_1}-x_{i_2})^{2m+1}$ is a generator of $Q_m(n,\k)$, as it is the lowest degree quasi-invariant in the sign module. However, from Lemma~\ref{lem:3.4} we have that in characteristic 3,
 \[(L+s_{23}L)K-(K+s_{23}K)L=c\prod_{i_1<i_2}(x_{i_1}-x_{i_2})^{2m+1},\]
 so $\prod_{i_1<i_2}(x_{i_1}-x_{i_2})^{2m+1}$ is not a generator. We can take this calculation further, and note that $(L+s_{23}L)K_1-(K+s_{23}K)L_1$ would then generate a copy of $\tign$, as the quotient of this module by the space generated by $(L+s_{23}L)K-(K+s_{23}K)L$ must be a trivial module.
 \end{Remark}

 It remains to consider the modules $\tign$, $\six$. To motivate the results that follow, we start by considering 0-quasi-invariants.

 \begin{Example}
 Note that from Corollary~\ref{cor:5gens} we know that $Q_0(3,\mathbf{F}_3)$ has 5 generators $1$, $x_1-x_2$, $(x_1-x_2)x_3$, $x_1$, $x_1(x_2+x_3)$ with no relations between them. By examining the dimension of the space of all homogeneous degree 3 polynomials, we have that $Q_0(3,\mathbf{F}_3)[3]$ is 10-dimensional. Since $\mathbf{F}_3[x_1,x_2,x_3]^{S_3}$ is 3-dimensional in degree $3$, 2-dimensional in degree $2$, and 1-dimensional in degree $1$, so far we have accounted for only $3+2+2+1+1=9$ dimensions. Moreover, every irreducible representation is accounted for, so this extra dimension must be an extension of an existing indecomposable representation. The only indecomposable representations that have nontrivial extensions are the $\triv$ generated by $x_1x_2x_3$ and the $\tign$ generated by
 \begin{align*}
 E :={}& (x_1x_2+x_1x_3+x_2x_3)x_1+(x_1+x_2+x_3)(x_1(x_2+x_3))\\
 ={}& {-}x_1^2x_2-x_1^2x_3+x_1x_2^2+x_1x_3^2.
 \end{align*}
 Indeed, the $\tign$ generated by $E$ extends to a $\six$ generated by
 \[F:=(x_1-x_2)x_1x_2.\]
 \end{Example}

 We will later see that the polynomials $E$, $F$ defined above are key to understanding $\tign$ and $\six$ in the quasi-invariants.

 \begin{Proposition}\label{prop:q0}
 $Q_0(3,\mathbf{F}_3)$ is freely generated by $1$, $x_1-x_2$, $(x_1-x_2)x_3$, $x_1$, $x_1(x_2+x_3)$, $F$ as~a~$\mathbf{F}_3[x_1,x_2,x_3]^{S_3}$-module.
 \end{Proposition}
 \begin{proof}
 We already know that the first 5 polynomials are independent. Now, let
 \[P_1+P_2(x_1-x_2)+P_3(x_1-x_2)x_3+P_4x_1+P_5(x_2+x_3)x_1+P_6F=0\]
 for symmetric polynomials $P_j$. Apply $1-s_{12}$ to this equation to get
 \[(P_4-P_2)(x_1-x_2)+(P_5-P_3)(x_1-x_2)x_3-P_6F=0.\]
 Next, apply $1+s_{23}$ to get
 \[(P_2-P_4)(x_1+x_2+x_3)+(P_5-P_3)(x_1x_2+x_1x_3+x_2x_3)+P_6E=0.\]
 Finally, note that as $E$ can be written in terms of symmetric polynomial multiples of $x_1$, ${(x_2+x_3)x_1}$, this equation would be a relation between the first 5 generators of $Q_0(3,\mathbf{F}_3)$. We~have seen this is impossible, so we have $P_6=0$, and hence all of the $P_j$ must be 0.

 Let $Q_0'$ be the submodule of $Q_0(3,\mathbf{F}_3)$ generated by these 6 polynomials. Then as the polynomials freely generate $Q_0'$ as a $\F_3[x_1,x_2,x_3]^{S_3}$-module, we have that the Hilbert series of $Q_0'$~is
 \[\mathcal{H}(Q_0')=\big(1+2t+2t^2+t^3\big)\mathcal{H}\big(\F_3[x_1,x_2,x_3]^{S_3}\big)=\frac{1+2t+2t^2+t^3}{(1-t)(1-t^2)(1-t^3)}=\frac{1}{(1-t)^3}\]
 by the fundamental theorem of symmetric polynomials. This is exactly the Hilbert series of $Q_0(3,\mathbf{F}_3)$, so $Q_0'=Q_0(3,\mathbf{F}_3)$ and there are no more generators of $Q_0(3,\mathbf{F}_3)$.
 \end{proof}

 Similar to how we only considered polynomials in the $(-1)$-eigenspace of $s_{12}$ for $\siv$, we only consider generators in the $(-1)$-eigenspace of $s_{12}$ for $\six$ and polynomials in the $1$-eigenspace of $s_{23}$ for $\tign$. Note that this is sufficient to describe the roles of $\six,\tign$, as both modules are generated by an element satisfying their respective constraints.

 \begin{Lemma}\label{lem:sixgen}
 \quad
 \begin{enumerate}\itemsep=0pt
\item[$1.$] Let $T\in Q_m(3,\mathbf{F}_3)$ generate a copy of $\tign$. Then $T$ is the sum of a symmetric polynomial multiple of $E\prod_{i_1<i_2}(x_{i_1}-x_{i_2})^{2m}$ and a symmetric polynomial. Conversely, any symmetric polynomial multiple of $E\prod_{i_1<i_2}(x_{i_1}-x_{i_2})^{2m}$ generates a copy of $\tign$ in $Q_m(3,\mathbf{F}_3)$.
\item[$2.$] Let $T_1\in Q_m(3,\mathbf{F}_3)$ generate a copy of $\six$. Then $T_1$ is the sum of a symmetric polynomial multiple of $F\prod_{i_1<i_2}(x_{i_1}-x_{i_2})^{2m}$ and a symmetric polynomial multiple of $\prod_{i_1<i_2}(x_{i_1}-x_{i_2})^{2m+1}$. Conversely, any symmetric polynomial multiple of $F\prod_{i_1<i_2}(x_{i_1}-x_{i_2})^{2m}$ generates a copy of $\six$ in $Q_m(3,\mathbf{F}_3)$.
 \end{enumerate}
 \end{Lemma}

 \begin{proof}
 1. We first prove the lemma for $m=0$. Consider some $T$ as above, and note that $(1-s_{12})T$ is contained in the sign representation, so by Proposition~\ref{prop:trivsign} we have $(1-s_{12})T=P(x_1-x_2)(x_1-x_3)(x_2-x_3)$ for some symmetric polynomial $P$. Then note that $PE$, $T$ generate two copies of $\tign$ with the same sign subrepresentation, so they generate a copy of
 \[(\tign\oplus\tign)/\sign\cong\tign\oplus\triv.\]
 So $T$ is the sum of $PE$ and a symmetric polynomial, as claimed.

 Now, consider general $m$. By the above we have that any $T$ must be of the form $T=PE+Q$ for symmetric polynomials $P$, $Q$. Then since $T$ is $m$-quasi-invariant, we have $(1-s_{12})T=P(x_1-x_2)(x_1-x_3)(x_2-x_3)$ is divisible by $(x_1-x_2)^{2m+1}$. So $P$ is divisible by $(x_1-x_2)^{2m}$, and it must also be divisible by $\prod_{i_1<i_2}(x_{i_1}-x_{i_2})^{2m}$ since it is symmetric.

 The converse is clear.

 2. This proof is similar to part~(1). For $m=0$, any $T_1$ must have that $(1+s_{23})T_1$ is in a~$\tign$ representation, so $(1+s_{23})T_1=PE$ for some $P\in\mathbf{F}_3[x_1,x_2,x_3]^{S_3}$. Then~$T_1$,~$PF$ generate a copy of
 \[(\six\oplus\six)/\tign\cong\six\oplus\sign,\]
 which implies the result for $m=0$. Then the extension to general $m$ is the same as in part (1). The converse is clear, as before.
 \end{proof}

 Finally, we can prove Theorem~\ref{thm:1} for $p=3$.

 \begin{Theorem}
 $Q_m(3,\mathbf{F}_3)$ is freely generated by $1$, the two generators $K$, $L$ from Theorem~$\ref{justRens}$, the two generators $K_1$, $L_1$ from Proposition~$\ref{prop:trix}$, and the generator $F\prod_{i_1<i_2}(x_{i_1}-x_{i_2})^{2m}$ from Lemma~$\ref{lem:sixgen}$.
 \end{Theorem}

 \begin{proof}
 Let us first show that there are no other generators of $Q_m(3,\mathbf{F}_3)$. Assume for contradiction that there is some other generator $T$ of $Q_m(3,\mathbf{F}_3)$. Then $T$ cannot generate a copy of $\triv$ by Proposition~\ref{prop:trivsign} and it cannot generate a copy of $\siv$ or $\trix$ by Theorem~\ref{justRens} and Proposition~\ref{prop:trix}. If it generates a copy of $\sign$, then by Proposition~\ref{prop:trivsign} it must be $\prod_{i_1<i_2}(x_{i_1}-x_{i_2})^{2m+1}$, but this polynomial is generated by $K$, $L$ by Lemma~\ref{lem:3.4}, so it cannot be a generator. If it generates a copy of $\tign$, then it is $E\prod_{i_1<i_2}(x_{i_1}-x_{i_2})^{2m}$ by Lemma~\ref{lem:sixgen}. But this is generated by $K_1$, $L_1$ by Remark \ref{rem:sign}. Finally, by Lemma~\ref{lem:sixgen} the only generator that generates a copy of $\six$ is $F\prod_{i_1<i_2}(x_{i_1}-x_{i_2})^{2m}$.

 Finally, we show there are no relations between the 6 generators. Note that this also implies $F\prod_{i_1<i_2}(x_{i_1}-x_{i_2})^{2m}$ is a generator, since it is not generated by the other 5 generators. But this is clear: we already know there are no relations between the first 5 generators by Corollary~\ref{cor:5gens}. If there was a relation involving $F\prod_{i_1<i_2}(x_{i_1}-x_{i_2})^{2m}$, then note that since every generator is generated by the generators of $Q_0(3,\mathbf{F}_3)=\mathbf{F}_3[x_1,x_2,x_3]$, this would induce a relation on those generators. Moreover, the generators other than $F\prod_{i_1<i_2}(x_{i_1}-x_{i_2})^{2m}$ each generate a copy of an indecomposable representation that is not $\six$, so they are each generated by the first 5 generators of $Q_0(3,\mathbf{F}_3)$. Meanwhile, $F\prod_{i_1<i_2}(x_{i_1}-x_{i_2})^{2m}$ is the only generator not generated by the first 5 generators, so the induced relation would be nontrivial. But there is no such relation by Proposition~\ref{prop:q0}.
 \end{proof}

 Note that these generators imply a Hilbert series that agrees with Theorem~\ref{thm:1} since $K$ is either a minimal degree Ren--Xu counterexample or has degree $3m+1$ if one does not exist. In~this~way, the Hilbert series of $Q_m(3,\mathbf{F}_3)$ agrees with that of $Q_m(3,\mathbf{Q})$ if and only if there does not exist a Ren--Xu counterexample. Ren--Xu counterexamples only exist when the conditions of Conjecture \ref{conj:ren} are satisfied, so Conjecture \ref{conj:ren} is also implied.

\subsection*{Acknowledgements}

We would like to thank the PRIMES program for making the project possible. We would also like to thank the referees for their useful comments and suggestions. This material is based upon work supported under the National Science Foundation Graduate Research Fellowship under Grant No.~2141064.

\pdfbookmark[1]{References}{ref}
\LastPageEnding


\begin{thebibliography}{99}
\footnotesize\itemsep=0pt

\bibitem{Alperin_1986}
Alperin J.L., Local representation theory: {M}odular representations as an
 introduction to the local representation theory of finite groups,
 \textit{Cambridge Stud. Adv. Math.}, Vol.~11,
 \href{https://doi.org/10.1017/CBO9780511623592}{Cambridge University Press},
 Cambridge, 1986.

\bibitem{ajs}
Andersen H.H., Jantzen J.C., Soergel W., Representations of quantum groups at a
 {$p$}th root of unity and of semisimple groups in characteristic {$p$}:
 independence of {$p$}, \textit{Ast\'erisque} \textbf{220} (1994), 323~pages.

\bibitem{beg}
Berest Yu., Etingof P., Ginzburg V., Cherednik algebras and differential
 operators on quasi-invariants,
 \href{https://doi.org/10.1215/S0012-7094-03-11824-4}{\textit{Duke Math.~J.}}
 \textbf{118} (2003), 279--337,
 \href{http://arxiv.org/abs/math.QA/0111005}{arXiv:math.QA/0111005}.

\bibitem{calogero}
Calogero F., Solution of the one-dimensional {$N$}-body problems with quadratic
 and/or inversely quadratic pair potentials,
 \href{https://doi.org/10.1063/1.1665604}{\textit{J.~Math. Phys.}} \textbf{12}
 (1971), 419--436.

\bibitem{CHALYKH2002193}
Chalykh O.A., Macdonald polynomials and algebraic integrability,
 \href{https://doi.org/10.1006/aima.2001.2033}{\textit{Adv. Math.}}
 \textbf{166} (2002), 193--259,
 \href{http://arxiv.org/abs/math.QA/0212313}{arXiv:math.QA/0212313}.

\bibitem{cmp/1104179957}
Chalykh O.A., Veselov A.P., Commutative rings of partial differential operators
 and {L}ie algebras, \href{https://doi.org/10.1007/BF02125702}{\textit{Comm.
 Math. Phys.}} \textbf{126} (1990), 597--611.

\bibitem{Etingof2000SymplecticRA}
Etingof P., Ginzburg V., Symplectic reflection algebras, {C}alogero--{M}oser
 space, and deformed {H}arish-{C}handra homomorphism,
 \href{https://doi.org/10.1007/s002220100171}{\textit{Invent. Math.}}
 \textbf{147} (2002), 243--348,
 \href{http://arxiv.org/abs/math.AG/0011114}{arXiv:math.AG/0011114}.

\bibitem{etingof2002lectures}
Etingof P., Strickland E., Lectures on quasi-invariants of {C}oxeter groups and
 the {C}herednik algebra, \textit{Enseign. Math.} \textbf{49} (2003), 35--65,
 \href{http://arxiv.org/abs/math.QA/0204104}{arXiv:math.QA/0204104}.

\bibitem{felder2001action}
Felder G., Veselov A.P., Action of {C}oxeter groups on {$m$}-harmonic
 polynomials and {K}nizhnik--{Z}amolodchikov equations,
 \href{https://doi.org/10.17323/1609-4514-2003-3-4-1269-1291}{\textit{Mosc.
 Math.~J.}} \textbf{3} (2003), 1269--1291,
 \href{http://arxiv.org/abs/math.QA/0108012}{arXiv:math.QA/0108012}.

\bibitem{MOSER1975197}
Moser J., Three integrable {H}amiltonian systems connected with isospectral
 deformations,
 \href{https://doi.org/10.1016/0001-8708(75)90151-6}{\textit{Adv. Math.}}
 \textbf{16} (1975), 197--220.

\bibitem{Ren_2020}
Ren M., Xu X., Quasi-invariants in characteristic {$p$} and twisted
 quasi-invariants,
 \href{https://doi.org/10.3842/SIGMA.2020.107}{\textit{SIGMA}} \textbf{16}
 (2020), 107, 13~pages,
 \href{http://arxiv.org/abs/1907.13417}{arXiv:1907.13417}.

\bibitem{ruij}
Ruijsenaars S.N.M., Complete integrability of relativistic {C}alogero--{M}oser
 systems and elliptic function identities,
 \href{https://doi.org/10.1007/BF01207363}{\textit{Comm. Math. Phys.}}
 \textbf{110} (1987), 191--213.

\bibitem{wang2022explicit}
Wang F., Toward explicit {H}ilbert series of quasi-invariant polynomials in
 characteristic {$p$} and {$q$}-deformed quasi-invariants, \textit{New York~J.
 Math.} \textbf{29} (2023), 613--634,
 \href{http://arxiv.org/abs/2201.06111}{arXiv:2201.06111}.

\end{thebibliography}
\end{document}